\documentclass{article}
\usepackage[a4paper, margin=2.5cm]{geometry}
\usepackage[title]{appendix}

\usepackage[dvipsnames]{xcolor}
\usepackage{hyperref}
\usepackage{amssymb}
\usepackage{amsmath}
\usepackage{enumerate}
\usepackage{threeparttable}
\usepackage{booktabs}
\usepackage{natbib}
  \setcitestyle{authoryear,open={(},close={)}}
\usepackage{tikz}
  \usetikzlibrary{positioning}
  \usetikzlibrary{shapes.geometric}
  \usetikzlibrary{fit}
  \usetikzlibrary{decorations.pathmorphing}
\usepackage{calc}

\usepackage{authblk}
\usepackage{subcaption}
\usepackage{cleveref}

\tikzset{
  nPick/.style={
    regular polygon,
    regular polygon sides=3,
    draw,
    fill = black,
    minimum size=3mm,
    inner sep=0pt,
    rotate=0
  },
  nDrop/.style={
    nPick,
    rotate=180
  },
  dMove/.style={
    draw,
    line width=1.6pt
  },
  dTravel/.style={
    draw,
    dashed,
    line width=1.2pt
  },
  dWait/.style={
    draw,
    dotted,
    line width=1.2pt
  },
  dTank/.style={
    draw=black!50!white
  },
  dProcess/.style={
    draw=black!25!white,
    line width=3.2pt
  },
  nTime/.style={
    draw=none,
    fill=none
  },
  nTank/.style={
    draw=none,
    fill=none,
    black!50!white
  }
}

\newcommand{\lpvar}[1]{\mathbf{#1}}
\newcommand{\bigM}{\operatorname{M}}
\title{A unified modeling framework and improved formulations for single-hoist cyclic scheduling}

\author{Mark\'o Horv\'ath\thanks{HUN-REN Institute for Computer Science and Control, Budapest, Hungary; marko.horvath@sztaki.hu; corresponding author}}

\begin{document}

\maketitle

\begin{abstract}
The cyclic hoist scheduling problem originates in electroplating lines, where a single or multiple hoists transport parts between processing tanks subject to technological constraints.
The objective is typically to determine a cyclic sequence of hoist movements that minimizes the cycle time while satisfying travel and processing constraints. 
Although the problem has been widely studied for several decades, the literature contains a puzzling phenomenon: different studies often report different optimal cycle times for the same benchmark instances, which limits the comparability and reproducibility of computational results.

In this paper, we revisit the modeling of cyclic hoist scheduling problems from a unified perspective. 
We introduce a consistent modeling approach for single-hoist problems and analyze several mixed-integer linear programming (MIP) formulations proposed in the literature.
Our analysis identifies modeling inconsistencies and clarifies the relationships between existing formulations. 
Based on these observations, we propose straightforward constraint programming (CP) models that can serve as baseline approaches, and we also derive improved MIP~formulations.
Extensive computational experiments compare the strength and performance of the investigated formulations.

To support reproducible research, we also provide a publicly available library containing benchmark instances and implementations of several CP and MIP~formulations for single-hoist cyclic hoist scheduling.
\end{abstract}

\section{Introduction}
The \emph{hoist scheduling problem} (HSP) originates in electroplating facilities, where hoists transfer parts between chemical baths according to their predefined processing sequences.
The problem also naturally arises in other material handling contexts, such as robotic cell scheduling (e.g., \citet{dawande2005sequencing, feng2024scheduling}).
Briefly stated, a single or multiple hoists transfer carriers containing the parts between tanks, and the scheduling of these movements must satisfy several technological constraints, including minimum and maximum processing times.
The \emph{cyclic hoist scheduling problem} (CHSP) aims to determine a repetitive sequence of hoist movements, called a \emph{cycle}, which is then executed continuously in order to process a large number of identical parts.
Determining a cycle with minimal duration is therefore directly related to maximizing the throughput of the production line.

50 years ago, \citet{phillips1976mathematical} published their paper in which the authors considered a problem to minimize the cycle time for a single hoist.
Since then, several solution approaches have been developed for the problem and its variants, including the multi-hoist case.
These approaches include
mixed-integer linear programming (e.g., \citet{li2014mixed, nait2016modeling, feng2018cyclic, mao2018mixed}),
constraint programming (e.g., \citet{wallace2020new, liu2025fast}),
other exact approaches (e.g., \citet{riera2002improved, che2004single, leung2006efficient, lei2014optimal}),
and heuristics or metaheuristic algorithms (e.g., \citet{chtourou2013hybrid, amraoui2016efficient, yan2016hybrid, laajili2021adapted, amraoui2024new}).

Despite the extensive literature, the examination of published results reveals a puzzling phenomenon.
Different studies that are supposed to solve the same basic problem often report different optimal cycle times for the same benchmark instances.
Table~\ref{tab:indicated} illustrates this issue for a basic single-hoist problem for several widely used instances, where the reported optimal values vary significantly across papers.
Even in the most recent articles, we can find results that are completely different from all previous ones.

\begin{table}
\centering
\caption{Reported values in the literature for a basic single-hoist problem on some benchmark instances.}
\label{tab:indicated}
\begin{threeparttable}[b]
\begin{tabular}{lrrrrrrrr}
\toprule
      & \multicolumn{8}{c}{Instances}\\
      \cmidrule{2-9}
Paper & PU & Mini & BO1 & BO2 & Cu & Zn & Ligne 1 & Ligne 2 \\
\midrule
\citet{phillips1976mathematical}    & 580 &   - &     - &     - &      - &      - &   - &   - \\
\citet{armstrong1994bounding} & 521 & - & 304.1 & 255.7 & 319.95 & 435.85 & - & - \\
\citet{leung2004optimal}            & 521 & 281 & 299.5 & 279.3 & 1847.2 & 1743.4 &   - &   - \\
\citet{che2011multi} & 521 & - & - & - & - & - & 393 & 712\\
\citet{che2015robust}               & 521 & 287 & 281.9 & 279.3 &      - &      - & 392 & 712 \\
\citet{mao2018mixed} & - & 284 & - & 279.3 & - & 1743.4 & - & - \\
\citet{laajili2021adapted}\tnote{*} & 521 & 340 & 333.2 & 279.3 & 1847.2 & 1743.4 & 425 & 712 \\
\citet{liu2025fast} & 518 & 284 & 267 & 257 & 1847 & 1742 & 392 & 722 \\
\bottomrule
\end{tabular}
\begin{tablenotes}
\item[-] not reported
\item[*] values are not indicated as optimal
\end{tablenotes}
\end{threeparttable}
\end{table}

Such discrepancies significantly limit the comparability of different approaches.
In many cases, new methods are evaluated primarily by comparing their results with previously reported cycle times.
If these reference values are inconsistent or incorrect, it becomes difficult to assess the true performance of new algorithms.
This issue is particularly problematic for heuristic methods, whose quality may be measured to incorrect reference values.

Importantly, these inconsistencies rarely arise due to fundamentally different problem definitions.
Hoist scheduling problems have been systematically classified by \citet{manier2003classification}, who introduced a unified notation, similar to the well-known classification scheme of \citet{graham1979optimization} for machine scheduling problems.
This notation allows researchers to precisely describe the main characteristics of the investigated hoist scheduling problem.
In principle, this should make it possible to define and reproduce the same problem settings across studies.

Although the notation used in the literature is relatively unified, the formalization of the problems is much less consistent.
In many studies, the constraints are described only in text, without providing an explicit mathematical formulation that clearly characterizes the set of feasible solutions.
To investigate the sources of these inconsistencies, we therefore focus primarily on constraint programming (CP) and mixed-integer linear programming (MIP) models.
In these approaches, the constraints must be explicitly specified in mathematical form, which makes it possible to clearly understand how the authors interpret the technological restrictions of the problem.
Thus, these formulations provide a transparent and reproducible basis for analyzing and comparing different models.
Our analysis of the literature suggests that several factors may contribute to the observed discrepancies.
In some cases, optimality claims appear to be incorrect; in others, published formulations contain modeling inaccuracies or omit important constraints.
Subtle differences in modeling assumptions may also lead to different feasible solution spaces, and benchmark instances themselves may be used in slightly different forms.
As a result, formulations that are intended to represent the same problem may in fact correspond to different variants.

These observations highlight the need for a systematic and reproducible modeling framework for cyclic hoist scheduling problems.
In this paper, we therefore revisit the modeling of cyclic hoist scheduling problems from a unified perspective.
Our main contributions are as follows: 
(1a)~a unified modeling approach for single-hoist cyclic scheduling;
(1b)~new constraint programming formulations;
(2a)~revision of existing mixed-integer linear programming formulations from the literature; 
(2b)~strengthened MIP~formulations; 
(3)~computational comparison of MIP~formulations; and 
(4)~a public library containing benchmark instances and solution approaches.

\paragraph{Unified modeling of single-hoist problems, new CP formulations.}
We propose a unified modeling approach for single-hoist cyclic scheduling problems that establishes a consistent notation and modeling structure.
This framework enables a systematic treatment of several important extensions, including multifunction tanks, multi-tank (multi-stage) operations, multi-degree problems (multi-cycles), and multi-part problems.
Compared to existing models in the literature, our approach provides a more thorough treatment of extensions such as multifunction tanks and multi-degree problems, and explicitly accounts for load–unload configurations, which are often neglected in previous formulations.
This framework also provides a foundation for future extensions to multi-hoist problems.
Our modeling observations naturally lead to straightforward CP formulations for both the simple and the multi-cycle variants of the problem, which can serve as baseline models.
For details, see \Cref{sec:shcs}.

\paragraph{Revised and strengthened MIP~formulations.}
We conduct a systematic revision of several MIP~formulations from the literature.
Our analysis identifies modeling inconsistencies and subtle differences in assumptions.
Building on this analysis, we derive strengthened MIP~formulations for the basic problem with simple cycles that combine the advantages of existing approaches.
For details, see \Cref{sec:forms}.

\paragraph{Computational comparison of MIP~formulations.}
We perform extensive computational experiments to evaluate and compare the strength and computational performance of the investigated formulations.
The results also clarify the optimal cycle times for several widely used benchmark instances.
For details, see \Cref{sec:comp}.

\paragraph{A public library for cyclic hoist scheduling.}
We publish a publicly available library\footnote{\url{https://github.com/hmarko89/CyclicHoistScheduling}} containing benchmark instances and implementations of several CP and MIP~formulations.
The aim of this library is to support reproducible and comparable research in cyclic hoist scheduling.
For details, see \Cref{sec:library}

\section{Hoist Scheduling}
In this section, we formally define and formulate the hoist scheduling problems considered in this paper. 
In \Cref{sec:hsp:main}, we provide a high-level overview to give a global picture, followed by a detailed description in \Cref{sec:hsp:prob}. 
\Cref{tab:hsp:notation} summarizes the notation used throughout the paper. 

While defining the problem, we draw on several works from the literature. 
In particular, papers that provide significant modeling contributions are discussed in \Cref{sec:hsp:prob}. 
For more comprehensive classifications and surveys, we refer the reader to \citet{manier2003classification, dawande2005sequencing, levner2010complexity, feng2017modelisation, feng2024scheduling} and \citet{laajili2021modelisation} (in French).

\subsection{Main overview}\label{sec:hsp:main}
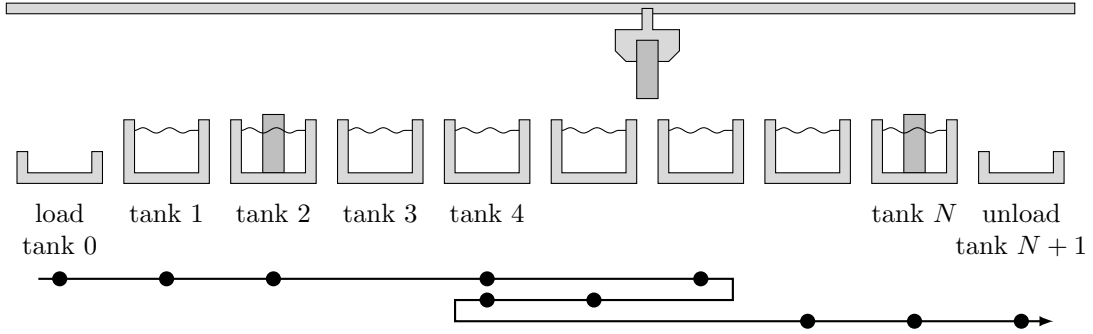
\begin{figure}
\centering
\begin{tikzpicture}
\begin{scope}[x=4,y=4]
\draw[fill=black!15!white] (0,5) rectangle (100,6);

\draw[fill=black!15!white] (59,0.5) --++(-1,0) --++(-1,1) --++(0,2) --++(2.5,0) --++(0,2) --++(1,0) --++(0,-2) --++(2.5,0) --++(0,-2) --++(-1,-1) --cycle;
\draw[fill=black!25!white] (59,2.5) --++(0,-5.5) --++(2,0) --++(0,5.5) --cycle;

\foreach \i in {2,8}{
\draw[fill=black!25!white] (10*\i+4,-4.5) --++(0,-5.5) --++(2,0) --++(0,5.5) --cycle;
}

\foreach \i in {1,...,8}{
\draw[fill=black!15!white] (10*\i+1,-5) --++(0,-6) --++(8,0) --++(0,6) --++(-1,0) --++(0,-5) --++(-6,0) --++(0,5) --cycle;
\draw[decorate,decoration={snake,amplitude=1}] (10*\i+2,-6) --++ (6,0);
}
\foreach \i in {0,9}{
\draw[fill=black!15!white] (10*\i+1,-8) --++(0,-3) --++(8,0) --++(0,3) --++(-1,0) --++(0,-2) --++(-6,0) --++(0,2) --cycle;
}
\node[below,align=center] at ( 5,-12) {load\\tank~$0$};
\node[below] at (15,-12) {tank~1};
\node[below] at (25,-12) {tank~2};
\node[below] at (35,-12) {tank~3};
\node[below] at (45,-12) {tank~4};
\node[below] at (85,-12) {tank~$N$};
\node[below,align=center] at (95,-12) {unload\\tank~$N+1$};

\draw[thick,->,-latex] (3,-20) -- (68,-20) -- (68,-22) -- (42,-22) -- (42,-24) -- (98,-24);
\foreach \i in {0,1,2,4,6}{ \node[circle,fill,minimum size=1mm,inner sep=2pt] at (10*\i+5,-20) {}; }
\foreach \i in {4,5}{ \node[circle,fill,minimum size=1mm,inner sep=2pt] at (10*\i+5,-22) {}; }
\foreach \i in {7,8,9}{ \node[circle,fill,minimum size=1mm,inner sep=2pt] at (10*\i+5,-24) {}; }

\end{scope}
\end{tikzpicture}
\caption{
A bi-directional line with a single hoist.
Currently, carriers are soaking in tanks~2 and~$N$, and another one is being transported by the hoist.
Tank~3 is not included in this processing sequence, however, the carrier re-enters multifunction tank~4.
}
\label{fig:hsp}
\end{figure}

Briefly stated, there is an electroplating line, consisting of a set of tanks, in which the parts are processed according to their predefined \emph{processing sequence}.
At the \emph{load station}, a part (or a batch of identical parts) is placed into a \emph{carrier} (e.g., a basket), which is then transported by a single or multiple \emph{hoists} between the tanks to soak according to the processing sequence, and is finally unloaded at the \emph{unload station}, see \Cref{fig:hsp}.
The hoists are identical and move on a single track. 
The hoist scheduling problem deals with the scheduling of the movements of the hoists.
In the cyclic hoist scheduling problem, the goal is to determine a repetitive sequence of hoist movements to process a larger quantity of identical parts.
The carriers are successively loaded onto the load station, and unloaded from the unload station.

In the rest of the paper, we consider cyclic problems.
For completeness, we mention that \citet{manier2003classification} distinguished 3~other problem classes, which are out of the scope of this paper.
\emph{Predictive} hoist scheduling problems comprise static problems in which a non-cyclic schedule is sought; for example, scheduling a part through the line to minimize the makespan (e.g., \citet{amraoui2016efficient}). 
In a \emph{dynamic} hoist scheduling problem, the schedule must be recomputed whenever a new part enters the line (e.g., \citet{feng2015dynamic, yan2017heuristic}).
In a \emph{reactive} hoist scheduling problem, decisions are made online by assigning tasks to hoists without a complete predefined schedule (e.g., \citet{jegou2006contract}).

\begin{table}
\centering
\caption{Notation for hoist scheduling problems.}
\label{tab:hsp:notation}
\begin{threeparttable}
\begin{tabular}{lcl}
\toprule
Parameter & Domain & Description\\
\midrule
$N$       & $\mathbb{Z}_{1\leq}$ & Number of chemical baths.\\
$n$       & $\mathbb{Z}_{1\leq}$ & Number of soaking operations.\\
$L_i$     & $\mathbb{R}_{0\leq}$ & Minimum duration for soaking operation~$i$.\\
$U_i$     & $\mathbb{R}_{0\leq}$ & Maximum duration for soaking operation~$i$, may be $\infty$.\\
$L_0$     & $\mathbb{R}_{0\leq}$ & Minimum duration to load a carrier to the load station.\\
$L_{n+1}$     & $\mathbb{R}_{0\leq}$ & Minimum duration to unload a carrier from the unload station.\\
$s_i$     & $\{0,1,\ldots,N+1\}$ & The tank of operation~$i$.\\
$S_i$     & $\subseteq\{1,\ldots,N\}$& The tanks of soaking operation~$i$ (multi-tank case).\\
$\operatorname{cap}(s_i)$ & $\mathbb{Z}_{1\leq}$ & Capacity of tank~$s_i$.\\
$d_i$     & $\mathbb{R}_{0\leq}$ & Duration of move~$i$.\\
$\tilde{e}_{i,j}$ & $\mathbb{R}_{0\leq}$ & Empty travel time from tank~$i$ to tank~$j$.\\
$e_{i,j}$ & $e_{i,j} = \tilde{e}_{s_i,s_j}$ & Empty travel time between the tanks of operation~$i$ and~$j$.\\
\bottomrule
\end{tabular}
\end{threeparttable}
\end{table}

\subsection{Problem definition}\label{sec:hsp:prob}
We consider a single line and index the tanks from~$0$ to~$N+1$ following the track.
That is, tank~$0$ is the load station, tanks~$1,\ldots,N$ are the chemical baths (or \emph{soaking tanks}), and tank~$N+1$ is the unload station, see \Cref{fig:hsp}.
In this case, the load and unload stations are \emph{dissociated}.
We consider a single \emph{part type} with a processing sequence consisting of $n+2$~operations, indexed from~$0$ to~$n+1$.
The first operation refers to loading the carrier to the load station $s_0 = 0$, and the last operation refers to unloading the carrier from the unload station $s_{n+1} = N+1$.
Operations $1,\ldots,n$ refer to soaking the carrier in their designated tanks~$s_1,\ldots,s_n$, respectively.
Multiple part types are discussed later.
Sometimes, load and unload stations are \emph{associated}, that is, loading and unloading are performed at the same station (i.e., $s_{n+1} = s_0$), however, in the single-hoist case, problems with associated load and unload stations can be converted into dissociated ones.
Each soaking tank is a \emph{monofunction tank}, that is, used to process only one operation of the processing sequence (i.e., $s_i \neq s_j$ for $1 \leq i < j \leq n$).
Only one carrier can soak in a tank at a time.
Multifunction tanks and multi-capacity tanks will be discussed later.

Operations must be processed without preemption.
The duration (or \emph{soak time}) of a soaking operation~$i$ ($1\leq i\leq n$) is bounded by minimum and maximum values, $L_i$ and $U_i$, respectively, which must be respected for quality reasons.
In some cases, the processing times are fixed, that is, the minimum and maximum durations coincide (i.e., $L_i = U_i$ for each soaking operation~$i$) (e.g., \citet{che2002multicyclic, leung2006efficient, che2010optimal}), which is typical in the case of robot scheduling problems (e.g. \citet{levner1997improved, yan2010branch}).
There might also be lower bounds~$L_0$ and~$L_{n+1}$ for loading and unloading, respectively, indicating the minimum required time to load/unload a carrier; however, upper bounds may have little practical use ($U_0 = U_{n+1} = \infty$).
Loading and unloading are not the responsibility of the hoists.
Loading and unloading might be performed simultaneously (which typically holds for dissociated load and unload stations) or not (which typically holds for associated load and unload stations).
Operations are subject to no-wait constraints; that is, once a part has been processed in a tank, it must be immediately transported to the next tank.

A hoist can only transport one carrier at a time.
In case of multiple hoists, the hoists are identical and move on the same track.
The empty travel time between tanks~$i$ and~$j$ is denoted by $\tilde{e}_{i,j}$.
The travel times are symmetric and satisfy the triangle inequality.
A \emph{move} is a loaded travel in which a hoist transports a carrier from one tank to the next.
A move actually consists of several sub-actions (e.g., lifting the carrier, dripping, transporting the carrier to the next tank, and lowering the carrier), however, the explicit modeling of these sub-actions is only necessary in the multi-hoist case to avoid collisions (e.g., \citet{che2013improved}).
The duration of move~$i$, i.e., the move associated with operation~$i$ ($0\leq i\leq n$), is denoted by~$d_i$.
Note that $\tilde{e}_{s_i,s_{i+1}} \leq d_i$.
Moves cannot be interrupted, and their durations are fixed.
That is, once a move starts (i.e., the corresponding operation ends), the carrier must be transported to the next tank without any intermediate pause due to technological reasons.
Very rarely, move durations are considered as lower bounds rather than exact values (e.g., \citet{liu2002cyclic, amraoui2008mixed}).

In the cyclic hoist scheduling problem, a repetitive sequence of hoist movements is sought, called a \emph{cycle}.
Based on the number of hoists, we have a \emph{single-hoist cyclic scheduling problem} (SHCSP) or a \emph{multi-hoist cyclic scheduling problem} (MHCSP).
In each cycle, a new carrier enters the line, which leaves the line possibly several cycles later.
Multi-cycles are discussed later.
The objective is typically to minimize the \emph{cycle time} (or \emph{period}), i.e., the length of the cycle.
Other objectives have also been considered (e.g., \citet{xu2004graph} minimized freshwater consumption),
sometimes as secondary objectives (e.g., \citet{che2015robust} maximized the robustness of the schedule, \citet{feng2014bi, yan2016hybrid} minimized the travel time of the hoist).
For the multi-hoist case, if the number of hoists is not fixed in advance, it may also be treated as an objective to minimize (e.g., \citet{leung2006efficient}).

\begin{figure}
\centering
\begin{tikzpicture}
\def\xmax{521}
\def\ymax{12}
\def\ctime{521}
\def\W{14cm}
\def\H{6cm}
\begin{scope}[x=\W/\xmax, y=\H/\ymax]
\foreach \h in {0,...,\ymax}{
  \draw[dTank] (0,\h) -- (\ctime,\h);
  \node[nTank, left] at (-1,\h) {\h};
}
\draw[dTank,dashed] (0,0) -- (0,\ymax);
\draw[dTank,dashed] (\ctime,0) -- (\ctime,\ymax);

\coordinate  (p0) at (  0, 0);
\coordinate  (d1) at ( 31, 1);
\coordinate (p10) at ( 41,10);
\coordinate (d11) at ( 68,11);
\coordinate  (r4) at ( 74, 4);
\coordinate  (p4) at ( 76, 4);
\coordinate  (d5) at (101, 5);
\coordinate  (p5) at (131, 5);
\coordinate  (d6) at (154, 6);
\coordinate (p11) at (168,11);
\coordinate (d12) at (190,12);
\coordinate  (p1) at (191, 1);
\coordinate  (d2) at (213, 2);
\coordinate (r12) at (216,12);
\coordinate (p12) at (220,12);
\coordinate  (d0) at (250, 0);
\coordinate  (p6) at (272, 6);
\coordinate  (d7) at (294, 7);
\coordinate  (p2) at (304, 2);
\coordinate  (d3) at (326, 3);
\coordinate  (r7) at (334, 7);
\coordinate  (p7) at (354, 7);
\coordinate  (d8) at (376, 8);
\coordinate  (p9) at (378, 9);
\coordinate (d10) at (425,10);
\coordinate  (p8) at (450, 8);
\coordinate  (d9) at (472, 9);
\coordinate  (p3) at (485, 3);
\coordinate  (d4) at (507, 4);
\coordinate  (r0) at (521, 0);
    
\foreach \h in {1,2,3}{
    \draw[dProcess,NavyBlue!25!white] (d\h) -- (p\h);
} 
\foreach \h in {4}{
    \draw[dProcess,NavyBlue!25!white] (d\h) -- (\ctime,\h);
    \draw[dProcess,BrickRed!25!white] (0,\h) -- (p\h);
}
\foreach \h in {5,6,7,8}{
    \draw[dProcess,BrickRed!25!white] (d\h) -- (p\h);
}
\foreach \h in {9}{
    \draw[dProcess,BrickRed!25!white] (d\h) -- (\ctime,\h);
    \draw[dProcess,OliveGreen!25!white] (0,\h) -- (p\h);
}
\foreach \h in {10}{
    \draw[dProcess,OliveGreen!25!white] (d\h) -- (\ctime,\h);
    \draw[dProcess,BurntOrange!25!white] (0,\h) -- (p\h);
}
\foreach \h in {11,12}{
    \draw[dProcess,BurntOrange!25!white] (d\h) -- (p\h);
}
\draw[thick,BurntOrange!25!white] (d0) -- (\ctime,0);

\draw[dMove,NavyBlue]    (p0) -- (d1);
\draw[dTravel]  (d1) -- (p10);
\draw[dMove,BurntOrange]   (p10) -- (d11);
\draw[dTravel] (d11) -- (r4);
\draw[dWait]    (r4) -- (p4);
\draw[dMove,BrickRed]    (p4) -- (d5);
\draw[dWait]    (d5) -- (p5);
\draw[dMove,BrickRed]    (p5) -- (d6);
\draw[dTravel]  (d6) -- (p11);
\draw[dMove,BurntOrange]   (p11) -- (d12);
\draw[dTravel] (d12) -- (p1);
\draw[dMove,NavyBlue]    (p1) -- (d2);
\draw[dTravel]  (d2) -- (r12);
\draw[dWait]   (r12) -- (p12);
\draw[dMove,BurntOrange]   (p12) -- (d0);
\draw[dTravel]  (d0) -- (p6);
\draw[dMove,BrickRed]    (p6) -- (d7);
\draw[dTravel]  (d7) -- (p2);
\draw[dMove,NavyBlue]    (p2) -- (d3);
\draw[dTravel]  (d3) -- (r7);
\draw[dWait]    (r7) -- (p7);
\draw[dMove,BrickRed]    (p7) -- (d8);
\draw[dTravel]  (d8) -- (p9);
\draw[dMove,OliveGreen]    (p9) -- (d10);
\draw[dTravel] (d10) -- (p8);
\draw[dMove,BrickRed]    (p8) -- (d9);
\draw[dTravel]  (d9) -- (p3);
\draw[dMove,NavyBlue]    (p3) -- (d4);
\draw[dTravel]  (d4) -- (r0);

\node[nPick,NavyBlue] at (p0) {};
\node[nDrop,NavyBlue] at (d1) {};
\node[nPick,NavyBlue] at (p1) {};
\node[nDrop,NavyBlue] at (d2) {};
\node[nPick,NavyBlue] at (p2) {};
\node[nDrop,NavyBlue] at (d3) {};
\node[nPick,NavyBlue] at (p3) {};
\node[nDrop,NavyBlue] at (d4) {};
\node[nPick,BrickRed] at (p4) {};
\node[nDrop,BrickRed] at (d5) {};
\node[nPick,BrickRed] at (p5) {};
\node[nDrop,BrickRed] at (d6) {};
\node[nPick,BrickRed] at (p6) {};
\node[nDrop,BrickRed] at (d7) {};
\node[nPick,BrickRed] at (p7) {};
\node[nDrop,BrickRed] at (d8) {};
\node[nPick,BrickRed] at (p8) {};
\node[nDrop,BrickRed] at (d9) {};
\node[nPick,OliveGreen] at (p9) {};
\node[nDrop,OliveGreen] at (d10) {};
\node[nPick,BurntOrange] at (p10) {};
\node[nDrop,BurntOrange] at (d11) {};
\node[nPick,BurntOrange] at (p11) {};
\node[nDrop,BurntOrange] at (d12) {};
\node[nPick,BurntOrange] at (p12) {};
\node[nDrop,BurntOrange] at (d0) {};

\node[nTime,below=2pt of p0] {0};
\node[nTime,above left=1pt and 1pt of d1] {31};
\node[nTime,above left=1pt and 1pt of p10] {41};
\node[nTime,above=1pt of d11] {68};
\node[nTime,below=1pt of p4] {76};
\node[nTime,above=1pt of d5] {101};
\node[nTime,below=1pt of p5] {131};
\node[nTime,below right=1pt and 1pt of d6] {154};
\node[nTime,above left=1pt and 1pt of p11] {168};
\node[nTime,above=1pt of d12] {190};
\node[nTime,below=1pt of p1] {191};
\node[nTime,below right=1pt and 1pt of d2] {213};
\node[nTime,above=1pt of p12] {220};
\node[nTime,below=2pt of d0] {250};
\node[nTime,above left=1pt and 1pt of p6] {272};
\node[nTime,above=1pt of d7] {294};
\node[nTime,below=1pt of p2] {304};
\node[nTime,below right=1pt and 1pt of d3] {326};
\node[nTime,below=1pt of p7] {354};
\node[nTime,below right=1pt and 1pt of d8] {376};
\node[nTime,above left=1pt and 1pt of p9] {378};
\node[nTime,above=1pt of d10] {425};
\node[nTime,below=1pt of p8] {450};
\node[nTime,above=1pt of d9] {472};
\node[nTime,below=1pt of p3] {485};
\node[nTime,above=1pt of d4] {507};
\node[nTime,below=2pt] at (\ctime,0) {\ctime};
\end{scope}
\begin{scope}[x=\W/\xmax, y=\H/\ymax, shift={(20,-2)}]
\draw[dMove] (0,0) -- (30,0);
\node[nPick] (f1) at (0,0) {};
\node[nDrop] at (30,0) {};
\node[nTime,right] at (40,0) {move};
\draw[dTravel] (120,0) -- (150,0);
\node[nTime,right] at (160,0) {empty travel};
\draw[dWait] (270,0) -- (300,0);
\node[nTime,right] at (310,0) {wait};
\draw[dProcess] (390,0) -- (420,0);
\node[nTime,right] (f2) at (430,0) {soaking};
\node[draw, rectangle,fit=(f1)(f2)] {};
\end{scope}
\end{tikzpicture}
\caption{
An optimal single-hoist schedule ($C = 521$) for the instance of \citet{phillips1976mathematical} with associated load and unload stations.
The colors refer to the four carriers used in the simple cycle.
}
\label{fig:pu:opt}
\end{figure}
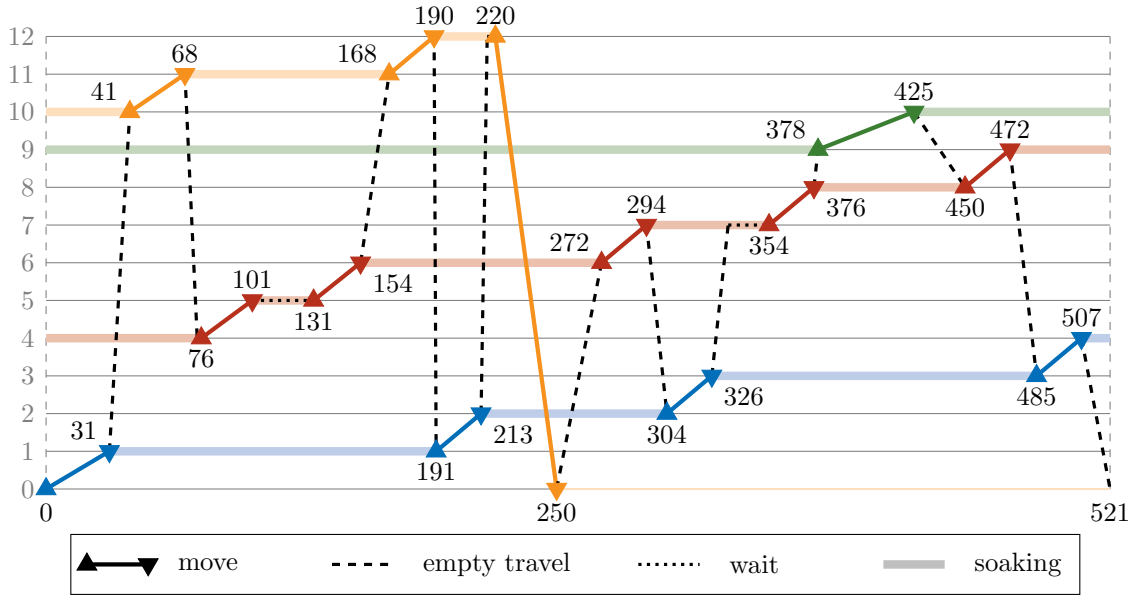

A solution to a cyclic hoist scheduling problem consists of three parts:
(i)~each move must be assigned to exactly one hoist;
(ii)~for each hoist, the execution order and timing of the assigned moves must be determined;
and (iii)~the movements between consecutive moves must be specified.
The solution can be illustrated by a \emph{time–way diagram}.
In the single-hoist case, the move-hoist assignment is trivial; and it is easy to see that, in addition to the loaded travel, two types of movements are sufficient: the \emph{empty travel} of the hoist to the tank of the next move; and idle \emph{waiting} at a tank.
\Cref{fig:pu:opt} indicates the time–way diagram of an optimal single-hoist solution for the instance proposed in~\citep{phillips1976mathematical}.
Note that in each cycle, four carriers are in use, and each carrier leaves the line three cycles after it enters.

In summary, the \emph{basic} single-hoist cyclic scheduling problem consists of a single line, and a simple cycle schedule is sought for a single hoist to process a single part type, where each operation is associated with a single-capacity monofunction tank.
Some possible extensions of the basic problem are discussed in the following.

\paragraph{Multifunction tanks.}
A soaking tank is a \emph{multifunction tank} (or \emph{reentrant tank}) if it is used to process multiple operations of the same processing sequence (i.e., $s_i = s_j$ for some $0 < i < j-1 < n$) (e.g., \citet{liu2002cyclic, feng2018cyclic}).
The operations associated with the same tank cannot overlap.
In some papers, using the same tank for different operations is called \emph{reentrance} or \emph{recirculation}.

\paragraph{Multi-degree problems, multi-part problems.}
In contrast to a \emph{simple cycle} (or \emph{1-degree cycle}), in an \emph{$r$-degree cycle} (or \emph{$r$-cycle}), $r$~carriers enter and leave the line in each cycle (e.g., \citet{che2011multi, li2014mixed, mao2018mixed}).

The multi-part case is analogous.
Multiple part types must be processed simultaneously on the same line, so that in each cycle, one carrier for each part type enters the line.
It is also possible that the different part types do not enter in a 1:1 ratio within a cycle, but rather according to a prescribed production ratio (e.g., \citet{zhao2013production}).
In a job shop environment, the different part types follow the same processing sequence (e.g., \citet{amraoui2008mixed, amraoui2012resolution, amraoui2013linear}); in a flow shop environment, each part type has its own processing sequence (e.g., \citet{feng2018cyclic}).

\paragraph{Multi-tank (multi-stage) operations, duplicated tanks.}
Operations with relatively long minimum durations may create bottlenecks. 
In such cases, \emph{duplicated tanks} (or \emph{parallel tanks}) are introduced to improve the throughput of the line (e.g., \citet{shapiro1988hoist}). 
That is, a \emph{multi-tank (multi-stage) operation}~$i$ is associated not with a single tank~$s_i$, but with a set~$S_i = \{s_i^1,\ldots,s_i^m\}$ of tanks, and each carrier must be soaked in exactly one of them. 
This setting is similar to the case where the operation is associated not with multiple single-capacity tanks, but with a single \emph{multi-capacity tank}, which allows multiple carriers to soak simultaneously.

Most often, multi-tank operations are handled within a simple cycle (e.g., \citet{zhou2003single}). 
In each cycle, one carrier enters the line. 
For a multi-tank operation~$i$, the carrier is assigned to an available tank in $S_i$, typically the one most recently vacated. 
The carrier is removed several cycles later, allowing it to soak over multiple cycles. 
Carriers are removed from the tanks in~$S_i$ in the same order in which they were inserted. 
In some cases, the number of tanks used for a multi-tank operation is not fixed but treated as a decision variable (e.g., \citet{liu2002cyclic, amraoui2012resolution, feng2018cyclic}).

\subsection{Related literature}
In the past decades, many variants of hoist scheduling problems have been investigated, for which several solution approaches have been developed.
For an overview, we refer the reader to the recent survey of \citet{feng2024scheduling}.
In the following, we review some papers that provide formulations for single-hoist problems and are thus relevant to our primary goal: to develop a unified modeling framework for single-hoist cyclic scheduling.
\Cref{tab:lit:mip} summarizes the main features of the reviewed MIP~formulations.

\begin{table}
\centering
\caption{Selected MIP~formulations for single-hoist scheduling problems from the literature.}
\label{tab:lit:mip}
\begin{threeparttable}
\begin{tabular}{lcccc}
\toprule
Paper & Degree & Load-Unload & Multifunction & Multi-tank \\
\midrule
\citet{phillips1976mathematical} & simple cycle & - & (yes) & -\\
\citet{liu2002cyclic} & simple cycle & associated & yes & variable \\
\citet{zhou2003single} & simple cycle & associated & (yes) & fixed \\
\citet{leung2004optimal} & simple cycle & - & - & -\\
\citet{feng2014bi} & simple cycle & - & - & - \\
\citet{che2015robust} & simple cycle & - & - & -\\
\citet{amraoui2012ptemporal} & simple cycle & associated & yes & variable \\
\citet{nait2016modeling} & simple cycle & associated & yes & variable \\

\addlinespace[1ex]
\citet{zhou2012mixed}     & multi-cycle   & - & - & - \\
\citet{li2014mixed}       & multi-cycle   & - & yes & - \\

\addlinespace[1ex]
\citet{amraoui2008mixed}  & 2-part     & - & - & - \\
\citet{amraoui2012resolution} & 2-part     & associated & yes & variable \\
\citet{amraoui2013linear} & multi-part & - & - & - \\
\citet{feng2018cyclic}    & multi-part & - & yes & variable \\
\bottomrule
\end{tabular}
\end{threeparttable}
\end{table}

\citet{phillips1976mathematical} proposed the first MIP~formulation to minimize the cycle time for a single part type in a simple cycle.
The authors also extended their model to multifunction tanks.
We revise the formulation in~\Cref{sec:shcsp:pu}, where we show that it is too restrictive; therefore, the resulting solution is not actually optimal in all cases.

\citet{shapiro1988hoist} developed a branch-and-bound solution procedure for the same problem setting.
Their model can be applied to lines with duplicated tanks, where all these tanks must be used in a solution, and can be generalized to determine multi-cycles. 

\citet{ng1995determining} considered the scenario where the number of duplicated tanks is not fixed.
The authors iteratively applied the solution method of \citet{shapiro1988hoist} to determine the optimal number of duplicated tanks.

\citet{liu2002cyclic} pointed out that the multifunction tank constraints of~\citet{phillips1976mathematical} do not prevent operations using the same tank from overlapping.
The authors proposed a corrected formulation for the problem, revised in~\Cref{sec:shcsp:liu}.
They considered a scenario where the move times are not exact, but lower bounds.
The authors also extended their formulation to multi-tank stages, where the number of tanks used for a multi-tank operation is a decision variable.

\citet{zhou2003single} extended the model of \citet{phillips1976mathematical} to multi-tank operations, where the number of tanks used for a multi-tank operation is fixed.
The formulation, revised in~\cref{sec:shcsp:zhou}, uses a correct calculation of the cycle time.
The corresponding multifunction tank constraints are not specified in that paper.

\citet{leung2004optimal} tightened the MIP~formulation of \citet{phillips1976mathematical}, and introduced a set of valid inequalities for the convex hull of the feasible solutions, see \Cref{sec:shcsp:leung}.
However, similarly to the original model, this formulation does not guarantee optimal solutions in certain cases.

\citet{amraoui2012ptemporal, nait2016modeling} proposed MIP~formulations for the basic problem as well as extensions to multifunction tanks and multi-tank operations.

\citet{feng2014bi} considered minimizing the travel time of the hoist as a secondary objective assuming non-Euclidean travel times.

\citet{che2015robust} adjusted the formulation of \citet{liu2002cyclic} to a bi-objective MIP~formulation for a robust optimization problem.
This model was then used to minimize the cycle time for a given level of robustness, where delays may occur in the execution of the moves or in the empty travels.

\citet{wallace2020new} proposed a non-linearized CP formulation capable of handling multiple tank capacities.

\citet{liu2025fast} emphasized the substantial potential of modern constraint programming solvers for hoist scheduling problems.
The authors transformed the classical single-hoist MIP~formulation (cf. \citet{zhou2003single}) into a CP~formulation, and solved the widely-used benchmark instances with OR-Tools CP-SAT.
However, based on the reported cycle times (see \Cref{tab:indicated}), there appear to be discrepancies: either the implementation is incorrect, or the benchmark dataset differs from those available in the literature.
For further discussion, we refer to \Cref{sec:comp:inst:bench}.

\citet{zhou2012mixed} developed the first MIP~formulation for processing a single part type in multi-cycles.
The model of \citet{li2014mixed} is applicable in the presence of multifunction tanks.
These two papers follow different modeling approaches, which we will discuss later in \Cref{sec:model:mc}.

\citet{amraoui2008mixed} proposed an MIP~formulation for two part types with the same processing sequence.
The model of \citet{amraoui2013linear} is applicable for an arbitrary number of part types that follow the same processing sequence.
\citet{amraoui2012resolution} provided extensions for the 2-part formulations, including multifunction tanks and multi-tank operations.

\citet{zhao2013production} considered a multi-part hoist scheduling problem in which the number of carriers entering the line per cycle for each part type is determined according to prescribed production ratios.

\citet{feng2018cyclic} developed an MIP~formulation for multiple part types with different processing sequences, which may contain multi-capacity multifunction tanks.

Some works are closely related to our primary goal, as they introduce systematic modeling approaches.
However, our approach for the single-hoist case covers a broader range of extensions and additionally provides reproducible formulations and benchmark instances to support systematic comparison.
\citet{amraoui2012ptemporal, nait2016modeling} used P-temporal Petri Nets to formally model and analyze single-hoist cyclic schedules for the basic problem, as well as for extensions involving multifunction tanks and multi-tank operations.
Recently, \citet{feng2024scheduling} provided a comprehensive review of robotic cell scheduling. 
In addition, the authors drew on several of these modeling approaches to present an MIP formulation for the basic problem and for extensions such as multifunction tanks, multi-tank operations, multi-cycle settings, and multiple hoists.

\section{Modeling single-hoist cyclic scheduling problems}\label{sec:shcs}
In this section, we address the modeling of single-hoist cyclic scheduling problems.
In~\Cref{sec:shcsp:pre}, we introduce the necessary definitions and notation.
\Cref{sec:model:sc,sec:model:mc} present the formulations for simple cycles and multi-cycles, respectively.
Before that, we describe the considered problem settings using two complementary classification schemes. 

First, we characterize the problem settings according to the classification of \citet{manier2003classification} (with the corresponding notation in parentheses).
We consider a single line ($\beta_1 = \beta_2 = \beta_3 = \emptyset$), a single hoist ($\beta_{41} = \emptyset$), and an arbitrary number of soaking tanks ($\beta_{42} = N$) with unit capacity ($\beta_{43} = \emptyset$).
We also consider an extension in which bottleneck operations are associated with multi-capacity tanks ($\beta_{43} = ct$).
The number of available carriers can be either unlimited or limited ($\beta_{5} \in \{ \emptyset,nc\}$), and a storage is available for them; thus, empty carriers do not remain on the line ($\beta_8 = \emptyset$).
Load and unload stations can be either associated or dissociated ($\beta_9 \in \{\emptyset,ass,diss\}$), when relevant.
An infinite number of identical parts are processed on the line ($\delta_1 = \emptyset$, $\delta_2 = \emptyset$), with an arbitrary number of soaking operations ($\delta_3 = n$), which may use the same tanks multiple times ($\delta_5 \in \{\emptyset,recr\}$).
After being unloaded, the carriers do not require cleaning ($\delta_4 = \emptyset$).
We also consider the multi-part case ($\delta_2 = nps$).
The objective is to find a cyclic schedule ($\alpha = CHSP$) with a minimum cycle time ($\gamma = Tmin$), however, our modeling approach is mostly independent of the objective function.

Next, we describe the problem settings using the more recent classification of \citet{feng2024scheduling}.
\textbf{Tank settings:} in addition to basic systems, we consider extended systems with reentrant tanks, parallel tanks, or both.
\textbf{Number of hoists:} a single hoist is used to transport the carriers.
\textbf{Part processing time:} each operation is associated with a processing time window, where the minimum and maximum durations may coincide (i.e., fixed processing times).
\textbf{Part processing routine:} in addition to the single-part case (single-part flow shop), we consider the multi-part extension, where different part types may follow different processing sequences (multiple part-type flow shop or job shop).
\textbf{Scheduling mode:} we consider both one-degree cyclic scheduling (i.e., simple cycles) and multi-degree cyclic scheduling.
In the latter case, either a single part type is processed in a multi-cycle, or multiple part types are processed within a single cycle.
\textbf{Optimization criteria:} our modeling approach is mostly independent of the objective function; however, in all considered formulations, the cycle time is minimized.

\subsection{Preliminaries}\label{sec:shcsp:pre}
\Cref{tab:shcsp:notation} summarizes the notation that is used throughout \Cref{sec:shcs,sec:forms}.
The continuous \emph{start time variable}~$\lpvar{t}_i$ denotes the start time of the move~$i$ ($0\leq i\leq n$), and the \emph{cycle time variable}~$\lpvar{C}$ denotes the cycle time.
The binary \emph{ordering variable}~$\lpvar{y}_{i,j}$ indicates whether move~$j$ starts after move~$i$ ($0\leq i < j\leq n$), i.e., $\lpvar{t}_i < \lpvar{t}_j$.
Some formulations may interpret variables~$\lpvar{y}_{i,j}$ with $i>j$ as well. 
Clearly, $\lpvar{y}_{j,i} = 1 - \lpvar{y}_{i,j}$ holds for all $i\neq j$.
We call a formulation a \emph{base formulation} if it uses only the variables $\lpvar{C}$, $\lpvar{t}_i$, and $\lpvar{y}_{i,j}$.
An \emph{extended formulation} uses additional variables as well.
The binary variable~$\lpvar{z}_i$ indicates whether the move~$i$ is the latest move ($0\leq i\leq n$), and the continuous variable~$\lpvar{t}_{\max}$ denotes the maximum start time among all moves (i.e., the start time of the latest one).

Without loss of generality, we assume that the cycle starts with the immediate execution of move~0, that is, $\lpvar{t}_0 = 0$.
This implies that $\lpvar{y}_{0,j} = 1$ ($1\leq j\leq n$) and $\lpvar{z}_0 = 0$.
Moreover, this means that no move crosses the cycle boundary, that is, $\lpvar{t}_i + d_i \leq \lpvar{C}$ holds for each move~$i$.

To simplify the notation, let $e_{i,j} = \tilde{e}_{s_i,s_j}$ denote the empty travel time from the tank of operation~$i$ to the tank of operation~$j$.
Recall that the travel times satisfy the triangle inequality, and that $e_{i,i+1} \leq d_i$ holds ($0\leq i\leq n$); both properties are used in the formulations.
In some constraints, we use an appropriately large constant~$\bigM$, which must be an upper bound on the cycle time.
For this purpose, we apply the upper bound~(\ref{eq:ub}) defined in the next subsection.

\begin{table}
\centering
\caption{Notation for the SHCSP formulations.}
\label{tab:shcsp:notation}
\begin{threeparttable}
\begin{tabular}{lcl}
\toprule
Variables & Domain & Description\\
\midrule
$\lpvar{t}_i$      & $\mathbb{R}_{0\leq}$ & Start time of move~$i$ ($0\leq i \leq n$).\\
$\lpvar{C}$ & $\mathbb{R}_{0\leq}$ & Cycle time.\\
$\lpvar{y}_{i,j}$ & $\{0,1\}$ & Indicates whether move~j starts after move~i ($0\leq i,j\leq n$, $i\neq j$).\\
$\lpvar{z}_i$ & $\{0,1\}$ & Indicates whether move~i is the latest move ($0\leq i\leq n$).\\
$\lpvar{t}_{\max}$ & $\mathbb{R}_{0\leq}$ & Maximum start time among all moves.\\
\bottomrule
\end{tabular}
\end{threeparttable}
\end{table}

As defined in \Cref{sec:hsp:prob}, we consider dissociative load and unload stations (i.e., $s_0 = 0$ and $s_{n+1} = N+1$).  
This is not a restriction: in the single-hoist case, any problem instance with associative load and unload stations (i.e., $s_0 = s_{n+1} = 0$) can be converted to an equivalent instance with dissociative stations without changing the set of feasible solutions or the optimal cycle time.  
Specifically, a fictive unload station, $N+1$, is introduced for the unloading operation such that $\tilde{e}_{0,N+1} = 0$ and $\tilde{e}_{0,j} = \tilde{e}_{N+1,j} = 0$ ($0 \leq j \leq N$), i.e., the travel times from/to the load station are identical to those from/to the unload station.

\subsubsection{Bounds on the cycle time}\label{sec:shcsp:bounds}
We present some lower and upper bounds on the cycle time for the case of minimizing the cycle time.
Note, however, that due to some loading/unloading configurations (see \Cref{sec:model:load}) or due to multi-tank operations (see \Cref{sec:model:multicap}), these bounds may need to be adjusted.

One can construct a primitive feasible solution for the problem in the following way.
The hoist executes the moves in the prescribed processing order.
When the carrier is placed in a tank, the hoist waits there until the corresponding operation is finished with its minimum processing duration.
After the execution of the last move, the hoist travels back to the load station.
This schedule provides an upper bound on the optimal cycle time, see~(\ref{eq:ub}), which will be used as the big-M constant or variable upper bounds in the formulations.
\begin{equation}\label{eq:ub}
\operatorname{UB} = d_0 + \sum_{i=1}^n (L_i + d_i) + e_{n+1,0}
\end{equation}

The hoist trajectory includes all moves and at least one empty travel.
Thus, the sum of move durations plus the minimum travel time to the load station is a lower bound on the cycle time.
For each soaking operation~$i$, the following are disjoint events: execution of move~$i-1$, processing of operation~$i$, execution of move~$i$, traveling (indirectly) from $s_{i+1}$ to $s_{i-1}$.
Thus, the total minimum duration of these events provides another lower bound on the cycle time.
These two observations can be combined in a single lower bound, see~(\ref{eq:lb}).
\begin{equation}\label{eq:lb}
\operatorname{LB} = \max\left\{ \min_{0\leq i\leq n}\{e_{i+1,0}\} + \sum_{i=0}^n d_i, \min_{1\leq i\leq n} \{ d_{i-1} + L_i + d_{i} + e_{i+1,i-1} \} \right\}
\end{equation}

\subsubsection{Time-way diagram}
The solution, i.e., the determined hoist movements, can be illustrated in a time-way diagram.
It is easy to see that in the single-hoist case, no other movements are needed in additional to the following ones:
(i)~the hoist executes a move;
(ii)~the hoist travels empty from the tank where the previous move is finished to the tank where the next move starts;
(iii)~the hoist waits either at the tank where the previous move is finished or at the tank where the next move starts.
Based on these observations, the start times of the moves and the cycle time are sufficient to describe a solution for a single-hoist cyclic scheduling problem.

For given move start times $t_0,\ldots,t_n$ and cycle time~$C$, the following procedure constructs the time-way diagram.

\paragraph{Step~1 (Operation start times).}
Let $t^+_0 = t_n + d_n$ be the start time of unloading.
Let $t^+_i = t_{i-1} + d_{i-1}\ (\operatorname{mod}\ C)$ be the start time of operation~$i$ ($1 \leq i \leq n$), that is, if $C \leq t_{i-1} + d_{i-1}$, then $t^+_i = t_{i-1} + d_{i-1} - C$.

\paragraph{Step~2 (Execution order).}
Let $\pi$ denote the execution order of moves, that is, $0 = t_{\pi_0} < t_{\pi_1} < \ldots < t_{\pi_n} < C$.

\paragraph{Step~3 (Trajectory).}
Following the execution order~$\pi$, we construct the parts of the hoist's trajectory.
For the incumbent move~$\pi_i$:
\begin{enumerate}[(1)]
\item Time interval $[t_{\pi_i},t_{\pi_i}+d_{\pi_i}]$: the hoist executes move~$\pi_i$.
\item Time interval $[t_{\pi_i}+d_{\pi_i},t_{\pi_i} + d_{\pi_i} + e_{\pi_i+1,\pi_{i+1}}]$ (may be empty): the hoist travels empty from tank~$s_{\pi_i+1}$ to tank $s_{\pi_{i+1}}$, where $s_{\pi_{n+1}}$ denotes the unload station.
Note that the time interval is empty if $s_{\pi_i+1} = s_{\pi_{i+1}}$.
\item Time interval $[t_{\pi_i} + d_{\pi_i} + e_{\pi_i+1,\pi_{i+1}},t^+_{\pi_{i+1}}]$ (may be empty): the hoist waits at tank $s_{\pi_{i+1}}$, where $t^+_{\pi_{n+1}} = C$.
\end{enumerate}

\subsection{Modeling simple cycles}\label{sec:model:sc}
A valid formulation for a single-hoist cyclic hoist scheduling problem with simple cycles must enforce four basic constraints:
\begin{itemize}
\item There must be enough time for empty travel between the consecutive moves, see \Cref{sec:model:travel}.
\item The cycle time must be long enough for the hoist to return to the load station after the latest move, see \Cref{sec:model:cycle}.
\item Minimum and maximum soak times must be respected, see \Cref{sec:model:soak}.
\item Loading and unloading times must be respected, see \Cref{sec:model:load}.
\end{itemize}
We also consider the following extensions:
\begin{itemize}
\item In case of multifunction tanks, operations using the same tank cannot overlap, see \Cref{sec:model:multi}.
\item Operations associated with multiple single-capacity tanks or a single multi-capacity tank, may be soaked for several cycles, see \Cref{sec:model:multicap}.
\item The number of carriers used in a cycle may be limited, see \Cref{sec:model:maxcarr}.
\end{itemize}

\subsubsection{Travel time constraints}\label{sec:model:travel}
There must be enough time between consecutive moves to allow the hoist to execute the former one and then travel empty to the next tank.
That is, if move~$j$ is the immediate successor of move~$i$, then $\lpvar{t}_i + d_i + e_{i+1,j}\leq \lpvar{t}_j$ must hold.
This can be enforced, for example, by the following indicator constraints~(\ref{eq:sc:travel:1})-(\ref{eq:sc:travel:0}).
\begin{align}    
\label{eq:sc:travel:1}\lpvar{y}_{i,j} = 1 &\to \lpvar{t}_j \geq \lpvar{t}_i + d_i + e_{i+1,j} & 0 \leq i < j \leq n\\
\label{eq:sc:travel:0}\lpvar{y}_{i,j} = 0 &\to \lpvar{t}_i \geq \lpvar{t}_j + d_j + e_{j+1,i} & 0 \leq i < j \leq n
\end{align}
These constraints imply the following lower bound on the start time variables: $\lpvar{t}_j \geq d_0 + e_{1,j}$ ($1\leq j\leq n$), since $\lpvar{t}_0 = 0$ and $\lpvar{y}_{0,j} = 1$.
Moreover, the actual interpretation of the $\lpvar{y}$-variables is the following: $\lpvar{y}_{i,j} = 1$ if $\lpvar{t}_i + d_i \leq \lpvar{t}_j$, and $\lpvar{y}_{i,j} = 0$ if $\lpvar{t}_j + d_j \leq \lpvar{t}_i$.

In an MIP~formulation, for example, the linearized constraint~(\ref{eq:travel:M}) can be used.
\begin{equation}
\label{eq:travel:M} \lpvar{t}_j \geq \lpvar{t}_i + d_i + e_{i+1,j} - \bigM(1-\lpvar{y}_{i,j}) \quad 0 \leq i,j \leq n:\ i\neq j
\end{equation}

\subsubsection{Cycle time constraints}\label{sec:model:cycle}
The cycle time must be long enough for the hoist to return to the load station after the latest move.
Constraint~(\ref{eq:sc:cycle}) ensures that this is fulfilled for every move.
\begin{equation}
\label{eq:sc:cycle}\lpvar{C} \geq \lpvar{t}_i + d_i + e_{i+1,0}\quad 0 \leq i \leq n
\end{equation}
Note that constraint~(\ref{eq:sc:cycle}) implies that each move starts and ends within the same cycle, that is, $\lpvar{t}_i + d_i \leq \lpvar{C}$ ($0 \leq i \leq n$).

\subsubsection{Soaking constraints}\label{sec:model:soak}

\begin{figure}
\centering
\begin{subfigure}[b]{0.45\textwidth}
\centering
\begin{tikzpicture}
\def\xmax{60}
\def\ymax{4}
\def\ctime{60}
\begin{scope}[x=6cm/\xmax, y=3cm/\ymax]
\foreach \h in {0,1,2,3,4}{
  \draw[dTank] (0,\h) -- (\xmax,\h);
}
\draw[dTank,dashed] (0,0) -- (0,\ymax);
\draw[dTank,dashed] (\ctime,0) -- (\ctime,\ymax);

\coordinate (p0) at (10,1);
\coordinate (d1) at (15,2);
\coordinate (p1) at (45,2);
\coordinate (d2) at (50,3);

\draw[dProcess] (d1) -- (p1);

\draw[dMove] (p0) -- (d1);
\draw[dMove] (p1) -- (d2);

\node[nPick] at (p0) {};
\node[nDrop] at (d1) {};
\node[nPick] at (p1) {};
\node[nDrop] at (d2) {};

\node[nTime,below=2pt of p0] {$\lpvar{t}_{i-1}$};
\node[nTime,above=2pt of d1] {$\lpvar{t}_{i-1} + d_{i-1}$};
\node[nTime,below=2pt of p1] {$\lpvar{t}_{i}$};

\node[nTime,below=2pt] at (0,0) {$0$};
\node[nTime,below=2pt] at (\ctime,0) {$\lpvar{C}$};
\end{scope}
\end{tikzpicture}
\caption{Soak time: $\lpvar{t}_{i} - (\lpvar{t}_{i-1} + d_{i-1})$.}
\label{fig:shcsp:soak:same}
\end{subfigure}
\hfill%
\begin{subfigure}[b]{0.45\textwidth}
\centering
\begin{tikzpicture}
\def\xmax{80}
\def\ymax{4}
\def\ctime{60}
\def\W{6cm}
\def\H{3cm}

\begin{scope}[x=\W/\xmax, y=\H/\ymax]
\foreach \h in {0,1,2,3,4}{
  \draw[dTank] (0,\h) -- (\xmax,\h);
  \node[nTank, left] at (-1,\h) {\h};
}
\draw[dTank,dashed] (0,0) -- (0,\ymax);
\draw[dTank,dashed] (\ctime,0) -- (\ctime,\ymax);

\coordinate (p0) at (25,1);
\coordinate (d1) at (30,2);
\coordinate (p1) at (10,2);
\coordinate (d2) at (15,3);
\coordinate (p1c) at (70,2);
\coordinate (d2c) at (75,3);

\draw[dProcess] (0,2) -- (p1);
\draw[dProcess] (d1) -- (p1c);

\draw[dMove] (p0) -- (d1);
\draw[dMove] (p1) -- (d2);
\draw[dMove] (p1c) -- (d2c);

\node[nPick] at (p0) {};
\node[nDrop] at (d1) {};
\node[nPick] at (p1) {};
\node[nDrop] at (d2) {};
\node[nPick] at (p1c) {};
\node[nDrop] at (d2c) {};

\node[nTime,below=2pt of p0] {$\lpvar{t}_{i-1}$};
\node[nTime,above=2pt of d1] {$\lpvar{t}_{i-1} + d_{i-1}$};
\node[nTime,below=2pt of p1] {$\lpvar{t}_{i}$};
\node[nTime,below=2pt of p1c] {$\lpvar{C} + \lpvar{t}_{i}$};

\node[nTime,below=2pt] at (0,0) {$0$};
\node[nTime,below=2pt] at (\ctime,0) {$\lpvar{C}$};
\end{scope}
\end{tikzpicture}
\caption{Soak time: $(\lpvar{C} + \lpvar{t}_{i}) - (\lpvar{t}_{i-1} + d_{i-1})$.}
\label{fig:shcsp:soak:diff}
\end{subfigure}
\caption{
Soaking either starts and ends within the same cycle ($\lpvar{y}_{i-1,i} = 1$, i.e., $\lpvar{t}_{i-1} < \lpvar{t}_i$) or it is in process at the beginning of the cycle ($\lpvar{y}_{i-1,i} = 0$, i.e., $\lpvar{t}_{i-1} > \lpvar{t}_i$).}
\label{fig:shcsp:soak}
\end{figure}
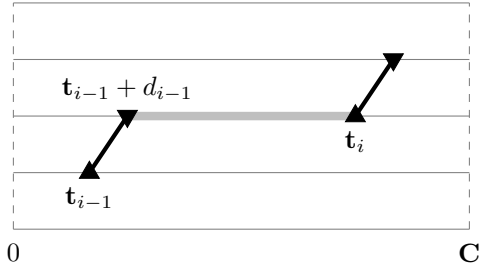
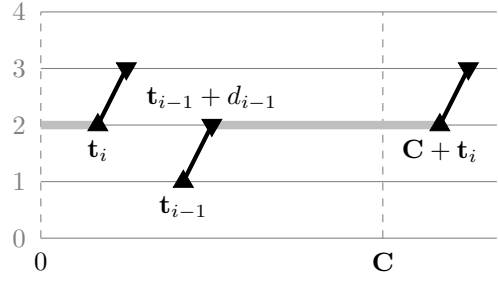

Minimum and maximum soak times must be respected.
For a soaking operation~$i$ ($1\leq i\leq n$), there are two possible scenarios:
\begin{enumerate}[(i)]
\item $\lpvar{t}_{i-1} + d_{i-1} \leq \lpvar{t}_i$: operation~$i$ starts and ends within the same cycle, see \Cref{fig:shcsp:soak:same}.
The soak time is $\lpvar{t}_i - (\lpvar{t}_{i-1} + d_{i-1})$.
\item $\lpvar{t}_i < \lpvar{t}_{i-1} + d_{i-1}$: operation~$i$ is in process at the beginning of the cycle, see \Cref{fig:shcsp:soak:diff}.
The soak time is $(\lpvar{C} + \lpvar{t}_i) - (\lpvar{t}_{i-1} + d_{i-1})$.
\end{enumerate}
Variable~$\lpvar{y}_{i-1,i}$ is suitable to distinguish these two cases for operation~$i$, see constraints~(\ref{eq:sc:soak:1})-(\ref{eq:sc:soak:0}).
\begin{align}
\label{eq:sc:soak:1}\lpvar{y}_{i-1,i} = 1 &\to L_i \leq \lpvar{t}_{i} - (\lpvar{t}_{i-1} + d_{i-1}) \leq U_i & 1 \leq i \leq n\\
\label{eq:sc:soak:0}\lpvar{y}_{i-1,i} = 0 &\to L_i \leq (\lpvar{C} + \lpvar{t}_{i}) - (\lpvar{t}_{i-1} + d_{i-1}) \leq U_i & 1 \leq i \leq n
\end{align}

\citet{chen1998cyclic} expressed these constraints with non-linear inequalities: $L_i \leq (\lpvar{y}_{i,i-1}\lpvar{C} + \lpvar{t}_i) - (\lpvar{t}_{i-1} + d_{i-1}) \leq U_i$ ($1\leq i\leq n$).
In an MIP~formulation, constraints~(\ref{eq:sc:soak:1}) and~(\ref{eq:sc:soak:0}) can be linearized, for example, using big-M formulations, see constraints~(\ref{eq:soak:M}) and~(\ref{eq:soakw:M}), respectively.
\begin{eqnarray}
\label{eq:soak:M} L_i - \bigM(1-\lpvar{y}_{i-1,i}) \leq \lpvar{t}_i - (\lpvar{t}_{i-1} + d_{i-1}) \leq U_i && 1 \leq i \leq n\\
\label{eq:soakw:M} L_i \leq (\lpvar{C} + \lpvar{t}_i) - (\lpvar{t}_{i-1} + d_{i-1}) \leq U_i + \bigM\lpvar{y}_{i-1,i} && 1 \leq i \leq n
\end{eqnarray}
As \citet{leung2004optimal} pointed out, no big-M term is needed for the upper bound inequality of constraint~(\ref{eq:soak:M}) and for the lower bound inequality of constraint~(\ref{eq:soakw:M}), since the corresponding inequalities hold in both of the two scenarios.
Our computational experiments in \Cref{sec:comp} will also demonstrate the strength of these improvements.

\subsubsection{Loading and unloading constraints}\label{sec:model:load}
Loading and unloading times must be respected.
We consider the following two load-unload configurations:
\begin{enumerate}[(i)]
\item Dissociated load and unload stations: the loading and unloading operations can be performed simultaneously.
\item Associated load and unload stations: the loading and unloading operations cannot be performed simultaneously.
Each carrier must be unloaded from the load–unload station before the next one can be loaded.
\end{enumerate}

In case of dissociated stations, the cycle time must be long enough to perform either operation:
\begin{align}
\label{eq:sc:load:diss}\max\{L_0,L_{n+1}\} \leq \lpvar{C}
\end{align}
Although such upper bounds may be of limited practical relevance, they can also be handled easily: $\lpvar{C} \leq \min\{U_0,U_{n+1}\}$.

In case of associated stations, there must be enough time for unloading the carrier and loading the next one:
\begin{align}
\label{eq:sc:load:ass} (\lpvar{t}_n + d_n) + L_{n+1} + L_0 \leq \lpvar{C} \end{align}
One can assume that $L_{n+1}=0$, since the unloading time can be merged into the loading time.
Note that in the presence of constraint~(\ref{eq:sc:load:ass}), the upper bound~(\ref{eq:ub}) must be adjusted by replacing the term $e_{n+1,0}$ with $\max\{e_{n+1,0},L_0+L_{n+1}\}$.
If upper bound~$U_0$ is also given, then constraint $\lpvar{C} \leq (\lpvar{t}_n + d_n) + U_0$ can be used.
However, if $e_{n+1,0} < U_0$, then the primitive schedule is not feasible, and thus the upper bound~(\ref{eq:ub}) is not valid.

\subsubsection{Multifunction tanks}\label{sec:model:multi}
The operations processed in the same multifunction tank cannot overlap.
Assume that operations~$i$ and~$j$ ($0< i< j-1< n$) require the same tank (i.e., $s_i = s_j$).
The former constraints ensure that an appropriate amount of time elapses between moves.
The remaining task is to ensure that the start and finish times of the operations are in the correct cyclic order.
There are four possible scenarios:
\begin{enumerate}[(i)]
\item $\lpvar{t}_{i-1} + d_{i-1} < \lpvar{t}_i < \lpvar{t}_{j-1} + d_{j-1} < \lpvar{t}_j$: both operations start and end within the cycle, operation~$i$ starts earlier;\label{enum:over:1}
\item $\lpvar{t}_j < \lpvar{t}_{i-1} + d_{i-1} < \lpvar{t}_i < \lpvar{t}_{j-1} + d_{j-1}$: operation~$j$ is in process at the beginning of the cycle;\label{enum:over:2}
\item $\lpvar{t}_{j-1} + d_{j-1} < \lpvar{t}_j < \lpvar{t}_{i-1} + d_{i-1} < \lpvar{t}_i$: both operations start and end within the cycle, operation~$j$ starts earlier;\label{enum:over:3}
\item $\lpvar{t}_i < \lpvar{t}_{j-1} + d_{j-1} < \lpvar{t}_j < \lpvar{t}_{i-1} + d_{i-1}$: operation~$i$ is in process at the beginning of the cycle\label{enum:over:4}
\end{enumerate}
These four cases correspond to the \emph{cyclic shifts} of the tuple $(\lpvar{t}_{i-1} + d_{i-1}, \lpvar{t}_i, \lpvar{t}_{j-1} + d_{j-1}, \lpvar{t}_j)$.
Therefore, exactly one of the four intervals $[\lpvar{t}_{i-1} + d_{i-1}, \lpvar{t}_i]$, $[\lpvar{t}_i, \lpvar{t}_{j-1} + d_{j-1}]$, $[\lpvar{t}_{j-1} + d_{j-1}, \lpvar{t}_j]$, $[\lpvar{t}_j, \lpvar{t}_{i-1} + d_{i-1}]$ is reversed (i.e., the corresponding event crosses the cycle boundary).
Thus, infeasible scenarios can be prevented, for example, with the following constraint (cf. \citet{feng2018cyclic}):
\begin{align}
\label{eq:sc:multi}\lpvar{y}_{i-1,i} + \lpvar{y}_{i,j-1} + \lpvar{y}_{j-1,j} + (1-\lpvar{y}_{i-1,j}) = 3
\end{align}
This implies $\lpvar{y}_{i-1,i} + \lpvar{y}_{j-1,j} \geq 1$, which means that operations cannot be simultaneously in process at the beginning of the cycle (i.e., carriers cannot soak simultaneously at the same tank at the beginning of the cycle).
Despite this, papers usually use two constraints, $\lpvar{y}_{i-1,i} + \lpvar{y}_{i,j-1} + \lpvar{y}_{j-1,j} + (1-\lpvar{y}_{i-1,j}) \geq 3$ and $\lpvar{y}_{i-1,i} + \lpvar{y}_{j-1,j} \geq 1$, to avoid such overlaps (e.g., \citet{liu2002cyclic,amraoui2012ptemporal}).

\begin{figure}
\centering
\begin{tikzpicture}
\def\xmax{290}
\def\ymax{3}
\def\ctime{290}
\def\W{12cm}
\def\H{3cm}

\begin{scope}[x=\W/\xmax, y=\H/\ymax]
\foreach \h in {0,...,\ymax}{
  \draw[dTank] (0,\h) -- (\xmax,\h);
  \node[nTank, left] at (-1,\h) {\h};
}
\draw[dTank,dashed] (0,0) -- (0,\ymax);
\draw[dTank,dashed] (\ctime,0) -- (\ctime,\ymax);

\coordinate (t0)   at (  0,0);
\coordinate (t0d0) at ( 10,1);
\coordinate (t1)   at ( 60,1);
\coordinate (t1d1) at ( 70,2);
\coordinate (t2)   at (180,2);
\coordinate (t2d2) at (190,1);
\coordinate (t3)   at (240,1);
\coordinate (t3d3) at (260,3);
\coordinate (t4)   at ( 80,3);
\coordinate (t4d4) at (100,1);
\coordinate (t5)   at (150,1);
\coordinate (t5d5) at (160,0);

\draw[dProcess,NavyBlue!25!white] (t0d0) -- (t1);
\draw[dProcess,NavyBlue!25!white] (t1d1) -- (t2);
\draw[dProcess,NavyBlue!25!white] (t2d2) -- (t3);
\draw[dProcess,NavyBlue!25!white] (t3d3) -- (\ctime,3);
\draw[dProcess,BrickRed!25!white] (0,3) -- (t4);
\draw[dProcess,BrickRed!25!white] (t4d4) -- (t5);
\draw[thick,BrickRed!25!white] (t5d5) -- (\ctime,0);

\draw[dWait] (t0d0) -- (t1);
\draw[dTravel] (t1d1) -- (t4);
\draw[dWait] (t4d4) -- (t5);
\draw[dTravel] (t5d5) -- (t2);
\draw[dWait] (t2d2) -- (t3);
\draw[dTravel] (t3d3) -- (\ctime,0);

\foreach \i in {0,1,2,3}{
  \draw[dMove,NavyBlue] (t\i) -- (t\i d\i);
  \node[nPick,NavyBlue] at (t\i) {};
  \node[nDrop,NavyBlue] at (t\i d\i) {};
}

\foreach \i in {4,5}{
  \draw[dMove,BrickRed] (t\i) -- (t\i d\i);
  \node[nPick,BrickRed] at (t\i) {};
  \node[nDrop,BrickRed] at (t\i d\i) {};
}

\node[nTime,below=2pt] at (0,0) {0};
\node[nTime,below=2pt] at (\ctime,0) {290};
\node[nTime,above=2pt of t0d0] {10};
\node[nTime,below=2pt of t1] {60};
\node[nTime,above left=1pt and 1pt of t1d1] {70};
\node[nTime,above=2pt of t4] {80};
\node[nTime,below=2pt of t4d4] {100};
\node[nTime,above=2pt of t5] {150};
\node[nTime,below=2pt of t5d5] {160};
\node[nTime,above=2pt of t2] {180};
\node[nTime,below=2pt of t2d2] {190};
\node[nTime,below=2pt of t3] {240};
\node[nTime,above=2pt of t3d3] {260};
\end{scope}
\end{tikzpicture}
\caption{
An example of a schedule with three operations that use the same tank.
A later operation (red carrier) appears between two earlier operations (blue carrier) in the cyclic timeline, 
showing that feasible schedules are not restricted to cyclic shifts of the original operation order.
}
\label{fig:ex2}
\end{figure}
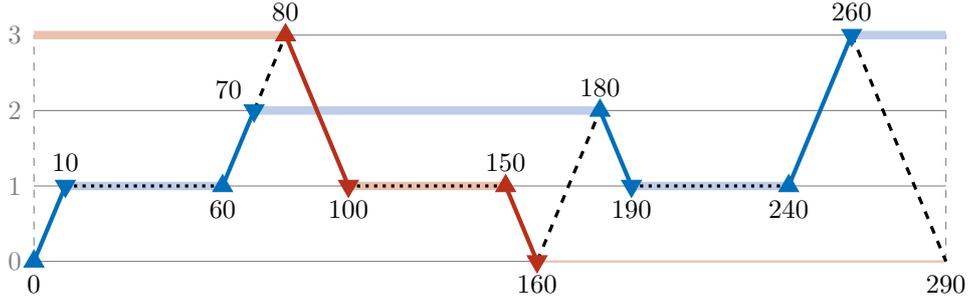

These pairwise no-overlap constraints can be also used if more than two soaking operations use the same tank; however, we show that this modeling idea cannot be generalized.
Assume that $s_{i_1} = s_{i_2} = \ldots = s_{i_r}$ for operations $i_1,i_2,\ldots,i_r$ ($i_p + 1 < i_{p+1}$ for each $1\leq p< r$).
Then, considering the cyclic shifts of the tuple $(\lpvar{t}_{i_1-1} + d_{i_1-1},\lpvar{t}_{i_1},\ldots,\lpvar{t}_{i_r-1} + d_{i_r-1},\lpvar{t}_{i_r})$ would lead to the following constraint:
\begin{align}\label{eq:sc:multi2}
\sum_{p=1}^r (\lpvar{y}_{i_p-1,i_p} + \lpvar{y}_{i_p,i_{p+1}-1}) = 2r-1
\end{align}
where the substitution $\lpvar{y}_{i_r,i_{r+1}-1} = (1-\lpvar{y}_{i_1-1,i_r})$ is applied to simplify the notation.
However, feasible schedules are not limited to cyclic shifts of this tuple, and therefore constraint~(\ref{eq:sc:multi2}) is not valid in general.
\Cref{fig:ex2} depicts an optimal schedule for an example, where $s = (0,1,2,1,3,1,0)$; $\tilde{e}_{i,j} = 10\times|i-j|$ ($0\leq i,j \leq 3$); $d_i = e_{i,i+1}$ ($0\leq i \leq 5$); $L_i = 50$ and $U_i = \infty$ ($1\leq i \leq 5$); $L_0 = L_{6} = 0$.
Notice that the operation sequence in multifunction tank~1 is the following: operation~1 (carrier~1), operation~5 (carrier~2), operation~3 (carrier~1),
which is not a cyclic shift of the operation order $(1,3,5)$.
This occurs because a later operation can be shifted across the cycle boundary, effectively being inserted between earlier operations in the cyclic timeline.

\subsubsection{Bottleneck operations: multi-capacity tanks, duplicated tanks}\label{sec:model:multicap}
\begin{figure}
\centering
\begin{tikzpicture}
\def\xmax{60}
\def\ymax{4}
\def\ctime{60}
\begin{scope}[x=12cm/\xmax, y=3cm/\ymax]
\foreach \h in {0,1,2,3,4}{
  \draw[dTank] (0,\h) -- (\xmax,\h);
}
\foreach \t in {0,10,20,30,40,50,60}{
  \draw[dTank,dashed] (\t,0) -- (\t,\ymax);
}

\coordinate (r1) at (5,1);
\coordinate (d1) at (8,1.7);
\coordinate (p1) at (32,1.7);
\coordinate (u1) at (35,3);
\coordinate (r2) at (15,1);
\coordinate (d2) at (18,2.0);
\coordinate (p2) at (42,2.0);
\coordinate (u2) at (45,3);
\coordinate (r3) at (25,1);
\coordinate (d3) at (28,2.3);
\coordinate (p3) at (52,2.3);
\coordinate (u3) at (55,3);
\coordinate (r4) at (35,1);
\coordinate (d4) at (38,1.7);
\coordinate (p4) at (2,1.7);
\coordinate (u4) at (5,3);
\coordinate (r5) at (45,1);
\coordinate (d5) at (48,2.0);
\coordinate (p5) at (12,2.0);
\coordinate (u5) at (15,3);
\coordinate (r6) at (55,1);
\coordinate (d6) at (58,2.3);
\coordinate (p6) at (22,2.3);
\coordinate (u6) at (25,3);

\foreach \i in {1,2,3}{
\draw[dProcess] (d\i) -- (p\i);
}
\foreach \i/\h in {4/1.7,5/2.0,6/2.3}{
\draw[dProcess] (d\i) -- (\ctime,\h);
\draw[dProcess] (0,\h) -- (p\i);
}
\foreach \i in {1,2,3,4,5,6}{
    \draw[dMove] (r\i) -- (d\i);
    \draw[dMove] (p\i) -- (u\i);
    \node[nPick] at (r\i) {};
    \node[nDrop] at (d\i) {};
    \node[nPick] at (p\i) {};
    \node[nDrop] at (u\i) {};
}
\node[nTank,left] at (0,1) {$s_{i-1}$};
\node[nTank,left] at (0,2) {$s_{i}$};
\node[nTank,left] at (0,3) {$s_{i+1}$};

\node[nTime,below=2pt of r1] {$\lpvar{t}_{i-1}$};
\node[nTime,below=2pt of p4] {$\lpvar{t}_i$};

\node[nTime,below=2pt] at (0,0) {$0$};
\node[nTime,below=2pt] at (10,0) {$\lpvar{C}$};
\node[nTime,below=2pt] at (20,0) {$2\lpvar{C}$};
\node[nTime,below=2pt] at (30,0) {$3\lpvar{C}$};
\node[nTime,below=2pt] at (40,0) {$4\lpvar{C}$};
\node[nTime,below=2pt] at (50,0) {$5\lpvar{C}$};
\node[nTime,below=2pt] at (60,0) {$6\lpvar{C}$};
\end{scope}
\end{tikzpicture}
\caption{
The process of operation~$i$ follows a 3-period pattern.
The carrier inserted into the tank in a given cycle is removed three cycles later, after being soaked for a duration between $2\lpvar{C}$ and $3\lpvar{C}$.
}
\label{fig:shcsp:multitank}
\end{figure}
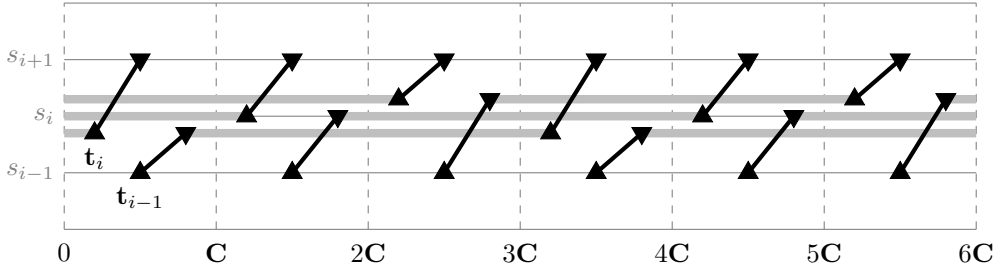

For this subsection, fix a \emph{bottleneck operation}~$i$ ($1\leq i\leq n$), which has a relatively long processing time.
Assume that $s_i$ is a monofunction tank with capacity greater than one, allowing a carrier to soak in the tank for multiple cycles.
The processing of the operation follows an $m$-period pattern for some fixed $1 \leq m \leq \operatorname{cap}(s_i)$, defined as follows.
In each cycle, a new carrier is inserted into the tank.
The carrier inserted into the tank in a given cycle is removed several cycles later, after being soaked for a duration between $(m-1)\lpvar{C}$ and $m\lpvar{C}$.
\Cref{fig:shcsp:multitank} depicts an example of an operation that is processed for three cycles.
Then, the corresponding constraints (\ref{eq:sc:soak:1})-(\ref{eq:sc:soak:0}) must be replaced with constraints (\ref{eq:msoak})-(\ref{eq:msoakw}).
\begin{align}
\label{eq:msoak}\lpvar{y}_{i-1,i} = 1 &\to L_i \leq ((m-1)\lpvar{C} + \lpvar{t}_{i}) - (\lpvar{t}_{i-1} + d_{i-1}) \leq U_i & \\
\label{eq:msoakw}\lpvar{y}_{i-1,i} = 0 &\to L_i \leq (m\lpvar{C} + \lpvar{t}_{i}) - (\lpvar{t}_{i-1} + d_{i-1}) \leq U_i &
\end{align}
Notice that variable $\lpvar{y}_{i-1,i}$ indicates whether the carrier that is inserted into the tank in a given cycle is removed $m-1$ or $m$ cycles later.
Also note that constraints (\ref{eq:sc:soak:1})-(\ref{eq:sc:soak:0}) are the same as constraints (\ref{eq:msoak})-(\ref{eq:msoakw}) for $m=1$.

In an MIP~formulation, linearized constraints (\ref{eq:msoak:M})-(\ref{eq:msoakw:M}) can be used.
\begin{eqnarray}
\label{eq:msoak:M} L_i - \bigM(1-\lpvar{y}_{i-1,i}) \leq ((m-1)\lpvar{C} + \lpvar{t}_i) - (\lpvar{t}_{i-1} + d_{i-1}) \leq U_i &&\\
\label{eq:msoakw:M} L_i \leq (m\lpvar{C} + \lpvar{t}_i) - (\lpvar{t}_{i-1} + d_{i-1}) \leq U_i + \bigM\lpvar{y}_{i-1,i} &&
\end{eqnarray}

\paragraph{Lower and upper bounds on the cycle time.}
Since the soak time of operation~$i$ is between $(m-1)\lpvar{C}$ and $m\lpvar{C}$, $(m-1)\lpvar{C} \leq U_i$ and $L_i \leq m\lpvar{C}$ must hold, otherwise the problem is infeasible.
This observation gives the following lower and upper bounds on the cycle time:
\begin{equation}\label{eq:multibounds}
\frac{L_i}{m} \leq \lpvar{C} \leq \frac{U_i}{m-1}
\end{equation}
The bounds~(\ref{eq:ub})-(\ref{eq:lb}) on the cycle time must be adjusted accordingly.

\paragraph{Duplicated tanks.}
If the bottleneck operation~$i$ is associated with multiple single-capacity tanks (duplicated tanks) instead of a single multi-capacity tank, the situation can be handled similarly.
The \emph{group} $S_i = \{s_i^1,s_i^2,\ldots,s_i^m\}$ of tanks where operation~$i$ can be processed can be merged into a single $m$-capacity tank, say~$s'_i$.
Each tank can belong to at most one group.

The advantage of this approach is that all carriers have identical soak times, ensuring consistent production quality.
However, move and travel times must be adjusted, which may lead to suboptimal solutions.
For instance, if $\tilde{e}_{s'_i,s_j} = \max_{1\leq k\leq m} \{\tilde{e}_{s_i^k,s_j}\}$ is used, i.e., the largest travel time among the tanks in the group, additional waiting times may be required when the hoist travels to or from tanks that are closer than the furthest one \citep{shapiro1988hoist}.

An alternative way of handling duplicated tanks will be discussed in \Cref{sec:model:mc}.

\paragraph{Determining the optimal number of duplicated tanks.}
The number~$m$ of duplicated tanks used for operation $i$ (i.e., the utilized portion of the capacity~$\operatorname{cap}(s_i)$) may also be treated as a decision variable.
Let the binary variable $\lpvar{u}_{i,m}$ indicate whether exactly $m$ units of the tank capacity are used ($1 \leq m \leq \operatorname{cap}(s_i)$).
Of course, exactly one value must be selected:
\begin{align}
\label{eq:mstage:sum}
\sum_{m=1}^{\operatorname{cap}(s_i)} \lpvar{u}_{i,m} = 1.
\end{align}
Now, the soak time of the operation also depends on the chosen capacity level:
\begin{align}
\label{eq:mstage:soak}
\lpvar{y}_{i-1,i} = 1 \text{ and } \lpvar{u}_{i,m} = 1 
&\to 
L_i \le ((m-1)\lpvar{C} + \lpvar{t}_i) - (\lpvar{t}_{i-1} + d_{i-1}) \le U_i,
&& 1 \le m \le \operatorname{cap}(s_i), \\
\label{eq:mstage:soakw}
\lpvar{y}_{i-1,i} = 0 \text{ and } \lpvar{u}_{i,m} = 1 
&\to 
L_i \le (m\lpvar{C} + \lpvar{t}_i) - (\lpvar{t}_{i-1} + d_{i-1}) \le U_i,
&& 1 \le m \le \operatorname{cap}(s_i).
\end{align}

For an MIP~formulation, the constraints can be linearized as:
\begin{align}
\label{eq:umsoak:M} L_i - \bigM(2-\lpvar{y}_{i-1,i} -\lpvar{u}_{i,m}) &\leq ((m-1)\lpvar{C} + \lpvar{t}_i) - (\lpvar{t}_{i-1} + d_{i-1}) \leq U_i + \bigM(1-\lpvar{u}_{i,m}) & 1\leq m\leq \operatorname{cap}(s_i)\\
\label{eq:umsoakw:M} L_i -\bigM(1-\lpvar{u}_{i,m}) &\leq (m\lpvar{C} + \lpvar{t}_i) - (\lpvar{t}_{i-1} + d_{i-1}) \leq U_i + \bigM(\lpvar{y}_{i-1,i} + 1 - \lpvar{u}_{i,m}) & 1\leq m\leq \operatorname{cap}(s_i)
\end{align}

Similar to the fixed-capacity scenario, move and travel times may require adjustment based on the number of tanks used \citep{liu2002cyclic}.

\subsubsection{Limited number of carriers}\label{sec:model:maxcarr}
The number of carriers available for production may be limited (e.g., \citet{wallace2020new}).
Constraint~(\ref{eq:maxcarr}) requires that at most $k$~carriers can be used in each cycle, i.e., exactly one carrier enters the line, and at most $k-1$ carriers soak at the beginning of the cycle.
\begin{equation}
1 + \label{eq:maxcarr}\sum_{i=1}^n (1-\lpvar{y}_{i-1,i}) \leq k
\end{equation}
Note that in case of multi-tank operations, constraint~(\ref{eq:maxcarr}) must be adjusted according to the tank capacities.
That is, for operation~$i$, the term $(1-\lpvar{y}_{i-1,i})$ must be replaced with $(m-\lpvar{y}_{i-1,i})$ or $(m\sum_{m=1}^{\operatorname{cap(s_i)}} \lpvar{u}_{i,m} - \lpvar{y}_{i-1,i})$, depending on whether the number of tanks/units, $m$, used for the operation is fixed or a decision variable.

\subsubsection{Valid inequalities}
As the $\lpvar{y}$-variables define a linear order, each formulation containing these variables may be strengthened by valid inequalities for the linear ordering polytope (see e.g., \citet{grotschel1985facets}).

\subsubsection{A straightforward CP formulation for simple cycles}\label{sec:sc:cp}
Our modeling considerations for simple cycles naturally lead to a straightforward CP~formulation.
The objective, for example, is to minimize the cycle time~$\lpvar{C}$, subject to constraints 
(\ref{eq:sc:travel:1})-(\ref{eq:sc:travel:0}), (\ref{eq:sc:cycle}), (\ref{eq:sc:soak:1})-(\ref{eq:sc:soak:0}), where $\lpvar{C} \in \{0,\ldots,\bigM\}$, $\lpvar{t}_0 = 0$, $\lpvar{t}_i \in \{0,\ldots,\bigM\}$ ($1\leq i\leq n$), $\lpvar{y}_{0,j} = 1$ ($1\leq j\leq n)$, $\lpvar{y}_{i,j} \in \{0,1\}$ ($1\leq i < j\leq n)$.
Depending on the load-unload configuration, constraint~(\ref{eq:sc:load:diss}) or (\ref{eq:sc:load:ass}) can be used.
In case of multifunction tanks, constraint~(\ref{eq:sc:multi}); in case of multi-tank operations, constraints~(\ref{eq:mstage:sum}) and (\ref{eq:mstage:soak})-(\ref{eq:mstage:soakw}) are added to the model.

\subsection{Modeling multi-degree cycles}\label{sec:model:mc}
In case of an $r$-cycle, $r$~carriers enter (and leave) the line, thus each operation~$i$ has $r$~\emph{copies}.
The modeling observations in \Cref{sec:model:sc} for simple cycles can be easily adapted for multi-degree problems.
In the following, we discuss only the modeling of the basic feasibility constraints and the extension to multifunction tanks.

Let operation~$(p,i)$ refer to the $p$th \emph{copy} of operation~$i$ ($1\leq p\leq r$, $0\leq i\leq n+1$).
That is, operation~$(p,i)$ refers to a carrier that enters the $p$th position in some cycle.
For moves $(p,i)$ and $(q,j)$, we use the lexicographic order: $(p,i) \leq (q,j)$ if and only if $p<q$, or $p=q$ and $i\leq j$.
Let the variable $\lpvar{t}_{(p,i)}$ denote the start time of move~$(p,i)$, and let the binary variable $\lpvar{y}_{(p,i),(q,j)}$ indicate whether $\lpvar{t}_{(p,i)} < \lpvar{t}_{(q,j)}$ ($(1,0)\leq (p,i) < (q,j) \leq (r,n)$).
Without loss of generality, we assume that the cycle starts with the execution of move~($1,0$), that is, $\lpvar{t}_{(1,0)} = 0$.
Consequently, $\lpvar{y}_{(1,0),(q,j)} = 1$ holds for each $(1,0)<(q,j)$.
Note that $d_{(p,i)} = d_i$, $e_{(p,i),(q,j)} = e_{i,j}$, $L_{(p,i)} = L_i$ and $U_{(p,i)} = U_i$, however, we use the general notation, because the formulation in the presented form will be easily adaptable to some extensions, see \Cref{sec:model:mc:adapt}.

\subsubsection{Symmetry breaking constraints}\label{sec:model:mc:symm}
Without loss of generality, we can assume that the carriers enter the line in the order $1,\ldots,r$:
\begin{align}
\label{eq:mc:symm}\lpvar{y}_{(p,0),(p+1,0)} = 1 && 1\leq p < r
\end{align}
However, we do not use this simplification in the following expressions, because the formulation in the presented form will be easily adaptable to the multi-part case, see \Cref{sec:model:mc:adapt}.

\subsubsection{Travel time constraints}
The travel time constraints are the following (cf. constraints~(\ref{eq:sc:travel:1})-(\ref{eq:sc:travel:0})):
\begin{align}
\label{eq:mc:travel1}\lpvar{y}_{(p,i),(q,j)} = 1 &\to \lpvar{t}_{(q,j)} \geq \lpvar{t}_{(p,i)} + d_{(p,i)} + e_{(p,i+1),(q,j)} & (1,0) \leq (p,i) < (q,j) \leq (r,n)\\
\label{eq:mc:travel0}\lpvar{y}_{(p,i),(q,j)} = 0 &\to \lpvar{t}_{(p,i)} \geq \lpvar{t}_{(q,j)} + d_{(q,j)} + e_{(q,j+1),(p,i)} & (1,0) \leq (p,i) < (q,j) \leq (r,n)
\end{align}

\subsubsection{Cycle time constraints}
The cycle time must be long enough for the hoist to return to the load station (cf. constraint~(\ref{eq:sc:cycle})):
\begin{align}
\label{eq:mc:cycle}\lpvar{C} &\geq \lpvar{t}_{(p,i)} + d_{(p,i)} + e_{(p,i+1),(p,0)} & (1,0)\leq (p,i)\leq (r,n)
\end{align}

\subsubsection{Soaking constraints}
The soaking constraints are the following (cf. constraints~(\ref{eq:sc:soak:1})-(\ref{eq:sc:soak:0})):
\begin{align}
\label{eq:mc:soak:1}\lpvar{y}_{(p,i-1),(p,i)} = 1 &\to L_{(p,i)} \leq \lpvar{t}_{(p,i)} - (\lpvar{t}_{(p,i-1)} + d_{(p,i-1)}) \leq U_{(p,i)} & (1,1) \leq (p,i) \leq (r,n)\\
\label{eq:mc:soak:0}\lpvar{y}_{(p,i-1),(p,i)} = 0 &\to L_{(p,i)} \leq (\lpvar{C} + \lpvar{t}_{(p,i)}) - (\lpvar{t}_{(p,i-1)} + d_{(p,i-1)}) \leq U_{(p,i)} & (1,1) \leq (p,i) \leq (r,n)
\end{align}
However, we must ensure that the different copies of an operation do not overlap.
The situation is very similar to the multifunction tank case discussed in \Cref{sec:model:multi}.
Therefore, the following pairwise constraints prohibit overlaps (cf. constraint~(\ref{eq:sc:multi})):
\begin{align}
\label{eq:mc:soak:over}\lpvar{y}_{(p,i-1),(p,i)} + \lpvar{y}_{(p,i),(q,i-1)} + \lpvar{y}_{(q,i-1),(q,i)} + (1-\lpvar{y}_{(p,i-1),(q,i)}) = 3 && (1,1) \leq (p,i) < (q,i) \leq (r,n)
\end{align}
This implies $\lpvar{y}_{(p,i-1),(p,i)} + \lpvar{y}_{(q,i-1),(q,i)} \geq 1$ for each pair $(p,q)$ of copies.
Consequently, $\sum_{p=1}^r \lpvar{y}_{(p,i-1),(p,i)} \geq r-1$ is a valid inequality.

Note that the no-wait constraints ensure that the entry order is preserved throughout the line, up to a cyclic shift; that is, $(1,2,\ldots,r)$, if the symmetry-breaking constraints provided in~\Cref{sec:model:mc:symm} are applied.
Consequently, a property that does not hold in the multifunction tank case, where operations may correspond to different processing stages, becomes valid here, as all considered operations belong to the same stage.
Specifically, feasible sequences for operation~$i$ correspond exactly to cyclic shifts of $((1,i),(2,i),\ldots,(r,i))$.
Therefore, constraint~(\ref{eq:mc:soak:over}) can be replaced by the following one (cf. constraint~(\ref{eq:sc:multi2})):
\begin{align}
\sum_{p=1}^r (\lpvar{y}_{(p,i-1),(p,i)} + \lpvar{y}_{(p,i),(p+1,i-1)}) = 2r-1 && 1\leq i\leq n
\end{align}
where the substitution $\lpvar{y}_{(r,i),(r+1,i-1)} = 1 - \lpvar{y}_{(1,i-1),(r,i)}$ is applied to simplify the notation.
Now, this constraint implies $\sum_{p=1}^r \lpvar{y}_{(p,i-1),(p,i)} \geq r-1$.
Despite this, these constraints are usually used together (e.g., \citet{zhou2012mixed, feng2024scheduling}).

\subsubsection{Loading and unloading constraints}
We consider the same load–unload configurations as in~\Cref{sec:model:load}.
Under the symmetry-breaking constraints from~\Cref{sec:model:mc:symm}, these constraints admit a simplified form, which we do not detail here.

If the load and unload stations are dissociated, constraints~(\ref{eq:mc:load:dl1})-(\ref{eq:mc:load:dl2w})
ensure that there is sufficient time between two first moves, and constraints~(\ref{eq:mc:load:du1})-(\ref{eq:mc:load:du2w}) enforce the same condition for last moves.
\begin{align}
\label{eq:mc:load:dl1}\lpvar{y}_{(p,0),(q,0)} = 1 \to\ & \lpvar{t}_{(p,0)} + L_{(q,0)} \leq \lpvar{t}_{(q,0)} & 1 \leq p < q \leq r\\
\label{eq:mc:load:dl1w}\lpvar{y}_{(p,0),(q,0)} = 1 \to\ & \lpvar{t}_{(q,0)} + L_{(p,0)} \leq \lpvar{C} + \lpvar{t}_{(p,0)} & 1 \leq p < q \leq r\\
\label{eq:mc:load:dl2}\lpvar{y}_{(p,0),(q,0)} = 0 \to\ & \lpvar{t}_{(q,0)} + L_{(p,0)} \leq \lpvar{t}_{(p,0)} & 1 \leq p < q \leq r\\
\label{eq:mc:load:dl2w}\lpvar{y}_{(p,0),(q,0)} = 0 \to\ & \lpvar{t}_{(p,0)} + L_{(q,0)} \leq \lpvar{C} + \lpvar{t}_{(q,0)} & 1 \leq p < q \leq r\\
\label{eq:mc:load:du1}\lpvar{y}_{(p,n),(q,n)} = 1 \to\ & \lpvar{t}_{(p,n)} + d_{(p,n)} + L_{(p,n+1)} \leq \lpvar{t}_{(q,n)} + d_{(q,n)} & 1 \leq p < q \leq r\\
\label{eq:mc:load:du1w}\lpvar{y}_{(p,n),(q,n)} = 1 \to\ & \lpvar{t}_{(q,n)} + d_{(q,n)} + L_{(q,n+1)} \leq \lpvar{C} + \lpvar{t}_{(p,n)} + d_{(p,n)} & 1 \leq p < q \leq r\\
\label{eq:mc:load:du2}\lpvar{y}_{(p,n),(q,n)} = 0 \to\ & \lpvar{t}_{(q,n)} + d_{(q,n)} + L_{(q,n+1)} \leq \lpvar{t}_{(p,n)} + d_{(p,n)} & 1 \leq p < q \leq r\\
\label{eq:mc:load:du2w}\lpvar{y}_{(p,n),(q,n)} = 0 \to\ & \lpvar{t}_{(p,n)} + d_{(p,n)} + L_{(p,n+1)} \leq \lpvar{C} + \lpvar{t}_{(q,n)} + d_{(q,n)} & 1 \leq p < q \leq r
\end{align}
Note that every second constraint accounts for cases where loading/unloading operations that cross the cycle boundary.

If the load and unload stations are associated, constraints~(\ref{eq:mc:load:a1})-(\ref{eq:mc:load:a4}) enforce a minimum gap between dropping off a carrier and picking up another one at the load-unload station.
\begin{align}
\label{eq:mc:load:a1}\lpvar{y}_{(p,n),(q,0)} = 1 \to\ & \lpvar{t}_{(p,n)} + d_{(p,n)} + L_{(p,n+1)} + L_{(q,0)} \leq \lpvar{t}_{(q,0)} & 1 \leq p < q \leq r\\
\label{eq:mc:load:a2}\lpvar{y}_{(p,n),(q,0)} = 0 \to\ & \lpvar{t}_{(p,n)} + d_{(p,n)} + L_{(p,n+1)} + L_{(q,0)} \leq \lpvar{C} + \lpvar{t}_{(q,0)} & 1 \leq p < q \leq r\\
\label{eq:mc:load:a3}\lpvar{y}_{(p,0),(q,n)} = 1 \to\ & \lpvar{t}_{(q,n)} + d_{(q,n)} + L_{(q,n+1)} + L_{(p,0)} \leq \lpvar{C} + \lpvar{t}_{(p,0)} & 1 < p \leq q \leq r\\
\label{eq:mc:load:a4}\lpvar{y}_{(p,0),(q,n)} = 0 \to\ & \lpvar{t}_{(q,n)} + d_{(q,n)} + L_{(q,n+1)} + L_{(p,0)} \leq \lpvar{t}_{(p,0)} & 1 < p \leq q \leq r
\end{align}

\subsubsection{Multifunction tanks}
Assume that operations~$i$ and~$j$ ($0< i< j-1< n$) are processed in the same tank (i.e., $s_i = s_j$).
Then, these operations of the same copy cannot overlap (cf. constraints~(\ref{eq:sc:multi})):
\begin{align}
\label{eq:mc:multi:same}\lpvar{y}_{(p,i-1),(p,i)} + \lpvar{y}_{(p,i),(p,j-1)} + \lpvar{y}_{(p,j-1),(p,j)} + (1-\lpvar{y}_{(p,i-1),(p,j)}) &= 3 & 1 \leq p\leq r
\end{align}
These operations of different copies also cannot overlap:
\begin{align}
\label{eq:mc:multi:diff1}\lpvar{y}_{(p,i-1),(p,i)} + \lpvar{y}_{(p,i),(q,j-1)} + \lpvar{y}_{(q,j-1),(q,j)} + (1-\lpvar{y}_{(p,i-1),(q,j)}) & = 3 & 1 \leq p < q \leq r\\
\label{eq:mc:multi:diff2}\lpvar{y}_{(p,j-1),(p,j)} + \lpvar{y}_{(p,j),(q,i-1)} + \lpvar{y}_{(q,i-1),(q,i)} + (1-\lpvar{y}_{(p,j-1),(q,i)}) & = 3 & 1 \leq p < q \leq r
\end{align}

\subsubsection{Adaptations to simple cycle extensions}\label{sec:model:mc:adapt}
Some extensions of the simple cycle problem can be addressed using multi-cycles.

\paragraph{Multi-tank operations.}
Consider a multi-tank operation~$i$ that can be performed in any tank from the group $S_i = \{s_i^1,\ldots,s_i^m\}$.
Rather than merging $S_i$ into a single $m$-capacity tank (cf. \Cref{sec:model:multicap}), the operation can be processed using an $m$-cycle, assigning the $k$th carrier to the $k$th tank in $S_i$.
This preserves the actual move and travel times ($d_{(p,i)}$ and $e_{(p,i),(q,j)}$).
However, without additional restrictions, the soak times of the different copies may vary, which may be undesirable.
If not all tanks must be used, an enumeration procedure may be needed to select the optimal number of duplicated tanks.
A major disadvantage of this approach is that handling multiple multi-tank operations can be inefficient, as the degree of the problem must be a common multiple of the group sizes.

\paragraph{Multiple part types.}
This model can be easily adapted to multiple part types, where a set of $r$~parts is processed simultaneously.
In each cycle, one carrier enters for each part type.
Here, $(p,i)$ refers to the $i$th operation of part type~$p$ ($1\leq p\leq r$, $0\leq i\leq n_r$), where $n_r$ denotes the number of soaking operations in the processing sequence of part type~$r$.
The symmetry-breaking constraints in \Cref{sec:model:mc:symm} no longer apply.

\subsubsection{Straightforward CP~formulations for multi-degree problems}
Our modeling observations for multi-cycles naturally lead to a CP~formulation for the single-part problem with multi-cycles.
The formulation can be easily adapted to the multi-part problem with a single cycle according to \Cref{sec:model:mc:adapt}.

The objective is, for example, to minimize the cycle time~$\lpvar{C}$
subject to constraints (\ref{eq:mc:symm}), (\ref{eq:mc:travel1})-(\ref{eq:mc:travel0}), (\ref{eq:mc:cycle}), and (\ref{eq:mc:soak:1})-(\ref{eq:mc:soak:over}),
where $\lpvar{C}\in\{0,\ldots,\bigM\}$,
$\lpvar{t}_{(1,0)}=0$, $\lpvar{t}_{(p,i)} \in\{0,\ldots,\bigM\}$ ($(1,0)< (p,i)\leq (r,n)$), $\lpvar{y}_{(1,0),(q,j)} = 1$ ($(1,0) < (q,j) \leq (r,n)$), $\lpvar{y}_{(p,i),(q,j)} \in \{0,1\}$ ($(1,0) < (p,i) \leq (q,j) \leq (r,n)$).
Depending on the load-unload configuration, constraints~(\ref{eq:mc:load:dl1})-(\ref{eq:mc:load:du2w}) or (\ref{eq:mc:load:a1})-(\ref{eq:mc:load:a4}) can be used.
In case of multifunction tanks, constraints (\ref{eq:mc:multi:same})-(\ref{eq:mc:multi:diff2}) are added to the model.
The constant~$\bigM$ can be chosen, for example, as $r$~times the upper bound~(\ref{eq:ub}) of the simple cycle.

\subsubsection{Other modeling approaches}
The MIP~formulation of \citet{li2014mixed} followed a different modeling approach.
In their formulation, move $(p,i)$ refers to the $k$th copy of operation~$i$ in chronological order within the cycle.
That is, $\lpvar{t}_{(1,i)} < \lpvar{t}_{(2,i)} < \ldots < \lpvar{t}_{(r,i)}$ holds for each operation~$i$ ($0\leq i\leq n$).
This ordering allowed certain constraints to be expressed more easily, but it was more cumbersome to generalize to multiple part types.
However, the proposed multifunction tank constraints are relatively complex.

\section{Revised and improved MIP~formulations for simple cycles}\label{sec:forms}
In this section, we revise several MIP~formulations from the literature for the basic problem.
That is, the problem setting considered consists of a single line with a single part type and a single hoist, an arbitrary number of single-capacity tanks, and an arbitrary number of operations, possibly with different load–unload configurations.
A simple cycle is sought with a minimal cycle time.

In \Cref{sec:shcsp:pu,sec:shcsp:leung,sec:shcsp:zhou}, we present the evolution of the extended formulations for the problem, including the models of \citet{phillips1976mathematical, leung2004optimal, zhou2003single}.
Combining their advantages, we present an improved formulation in \Cref{sec:shcsp:imp1}.
In \Cref{sec:shcsp:liu}, we revise the base formulation of \citet{liu2002cyclic}, and we propose a strengthened version in \Cref{sec:shcsp:imp2}.
The comparison of the investigated formulations is presented in~\Cref{sec:comp}.

\subsection{Formulation of \citet{phillips1976mathematical}}\label{sec:shcsp:pu}
\begin{subequations}
\renewcommand{\theequation}{A\arabic{equation}}
\citet{phillips1976mathematical} aimed to minimize the cycle time, for which they used three sets of constraints.
Although the authors did not introduce the cycle time variable~$\lpvar{C}$, we will use it here to simplify the notation (see constraint~(\ref{eq:pu:cvar})).
We emphasize that all our modifications leave the formulation equivalent to the original version.
%In \Cref{sec:shcsp:pu:multi}, we present how the authors extended their formulation to multifunction tanks.
%In \Cref{sec:shcsp:pu:remark}, we address our concerns with the formulation.

The objective is to minimize the cycle time~$\lpvar{C}$, where the cycle time is considered to be the return time of the hoist after the execution of the latest move, see constraint~(\ref{eq:pu:cvar}).
\begin{eqnarray}
\label{eq:pu:cvar} \lpvar{C} = \lpvar{t}_{\max} + \sum_{i=1}^n ( d_i + e_{i+1,0} ) \lpvar{z}_i
\end{eqnarray}
The first set of constraints, (\ref{eq:pu:3})-(\ref{eq:pu:2}), ensures that there is exactly one latest move, and that the variable~$\lpvar{t}_{\max}$ is properly defined.
\begin{eqnarray}
\label{eq:pu:3} \sum_{i=1}^n \lpvar{z}_i = 1 && \\
\label{eq:pu:1} \lpvar{t}_{\max} \geq \lpvar{t}_i && 1 \leq i \leq n\\
\label{eq:pu:2} \lpvar{t}_{\max} \leq \lpvar{t}_i + \bigM(1 - \lpvar{z}_i) && 1 \leq i \leq n
\end{eqnarray}
The second set of constraints, (\ref{eq:shscp:pu:4})-(\ref{eq:shscp:pu:6}), consists of the travel time constraints (cf. constraints~(\ref{eq:sc:travel:1})-(\ref{eq:sc:travel:0})).
Note that the case $j=i+1$ for constraint~(\ref{eq:shscp:pu:6}) is covered by the first inequality of constraint~(\ref{eq:pu:57}).
\begin{eqnarray}
\label{eq:shscp:pu:4}\lpvar{t}_j \geq (\lpvar{t}_i + d_i) + e_{i+1,j} - \bigM(1-\lpvar{y}_{i,j}) && 0 \leq i < j \leq n\\
\label{eq:shscp:pu:6}\lpvar{t}_i \geq (\lpvar{t}_j + d_j) + e_{j+1,i} - \bigM\lpvar{y}_{i,j} && 1 \leq i < j \leq n:\ i + 1 < j
\end{eqnarray}
The third set of constraints, (\ref{eq:pu:57})-(\ref{eq:pu:89}), consists of the soaking constraints (cf. constraints~(\ref{eq:soak:M})-(\ref{eq:soakw:M})).
\begin{eqnarray}
\label{eq:pu:57}L_i - \bigM(1-\lpvar{y}_{i-1,i}) \leq \lpvar{t}_i - (\lpvar{t}_{i-1} + d_{i-1}) \leq U_i + \bigM(1-\lpvar{y}_{i-1,i}) && 1 \leq i \leq n\\
\label{eq:pu:89}
L_i - \bigM\lpvar{y}_{i-1,i} \leq (\lpvar{C} + \lpvar{t}_i) - (\lpvar{t}_{i-1} + d_{i-1}) \leq U_i + \bigM\lpvar{y}_{i-1,i} && 1 \leq i \leq n
\end{eqnarray}

In summary, the formulation of \citet{phillips1976mathematical} aims to minimize the cycle time~$\lpvar{C}$, subject to constraints (\ref{eq:pu:cvar})-(\ref{eq:pu:89}), where $\lpvar{t}_0 = 0$, $0\leq \lpvar{t}_i$ ($1\leq i\leq n$), $\lpvar{y}_{0,j} = 1$ ($1\leq j\leq n$), $\lpvar{y}_{i,j} \in \{0,1\}$ ($1\leq i < j \leq n$), $\lpvar{z}_i \in \{0,1\}$ ($1\leq i \leq n$).

\subsubsection{Extension to multifunction tanks}\label{sec:shcsp:pu:multi}
The authors extended their formulation to multifunction tanks, by adding constraint~(\ref{eq:pu:overc}) to the model, and by replacing the appropriate constraints~(\ref{eq:shscp:pu:4})-(\ref{eq:shscp:pu:6}) with the following ones:
\begin{eqnarray}
\label{eq:pu:overc} \lpvar{y}_{i-1,i} + \lpvar{y}_{j-1,j} \geq 1 && 0 \leq i + 1 < j \leq n:\ s_i = s_j\\
\label{eq:pu:overij}\lpvar{t}_{j-1} \geq \lpvar{t}_{i-1} + d_{i-1} + L_i + d_i + e_{i+1,j} - \bigM(1-\lpvar{y}_{i-1,j-1})&& 0 \leq i + 1 < j \leq n:\ s_i = s_j\\
\label{eq:pu:overji}\lpvar{t}_{i-1} \geq \lpvar{t}_{j-1} + d_{j-1} + L_j + d_j + e_{j+1,i} - \bigM\lpvar{y}_{i-1,j-1}&& 0 \leq i + 1 < j \leq n:\ s_i = s_j
\end{eqnarray}
Constraint~(\ref{eq:pu:overc}) ensures that at least one of the operations is not in process at the beginning of the cycle.
Constraints~(\ref{eq:pu:overij})-(\ref{eq:pu:overji}) enforce a minimum gap between the operation start times for the corresponding operation end time.

\end{subequations}

\subsubsection{Remarks, issues}\label{sec:shcsp:pu:remark}
In the following, we discuss our observations regarding the formulation.

\paragraph{Unspecified loading and unloading constraints.}
Although the authors provided an example in which a lower bound was given for operation~0, their formulation did not include any constraints to enforce this.
This inconsistency can be easily fixed by adding constraint~(\ref{eq:sc:load:diss}) or~(\ref{eq:sc:load:ass}) to the model, depending on the load-unload configuration.

\paragraph{Narrowed solution space.}
To the best of our knowledge, this issue has not been raised before.
The authors fixed not only the beginning (recall that $\lpvar{t}_0 = 0$), but also the end of the cycle.
That is, they considered the expression $\lpvar{t}_{\max} + \sum_{i=1}^n \left( d_i + e_{i+1,0} \right) \lpvar{z}_i$ as the cycle time, which is the start time of the latest move, plus its execution time, plus the empty travel time to tank~0.
Thus, the cycle immediately ends with the return of the hoist to the load station after the execution of the latest move.
We call such a formulation \emph{restricted} with respect to the cycle finish.
This eliminates the possibility of the hoist waiting at the load station at the end of the cycle to further soak the parts in the tanks.
It has several consequences.
Clearly, this restricted structure narrows the set of feasible solutions, and therefore the resulting cycle time may be higher than the true optimum.
Moreover, if $e_{n+1,0} < L_0$ and constraint~(\ref{eq:sc:load:ass}) is imposed, then the primitive solution constructed in \Cref{sec:shcsp:bounds} may no longer be feasible. 
Consequently, the upper bound~(\ref{eq:ub}) may become invalid, and the formulation itself may even be infeasible.

\begin{table}
\centering
\caption{Instance \texttt{ex1}.}
\label{tab:cntex1}
\begin{tabular}{c|cccc|cccc}
\toprule
$i$ & $e_{i,0}$ & $e_{i,1}$ & $e_{i,2}$ & $e_{i,3}$ & $s_i$ & $d_i$ & $L_i$ & $U_i$\\
\midrule
0 &  0 & 10 & 20 &  0 & 0 & 10 &   0 & $\infty$\\
1 & 10 &  0 & 10 & 10 & 1 & 10 &  40 & 100\\
2 & 20 & 10 &  0 & 20 & 2 & 20 & 120 & $\infty$\\
3 &  0 & 10 & 20 &  0 & 0 &  - &   0 & $\infty$\\
\bottomrule
\end{tabular}
\end{table}

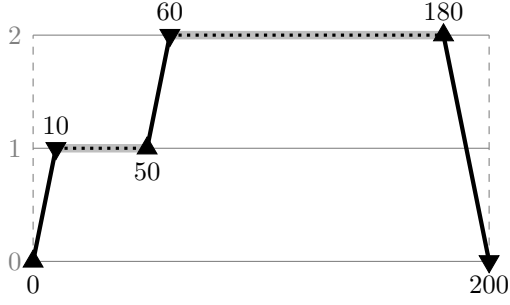
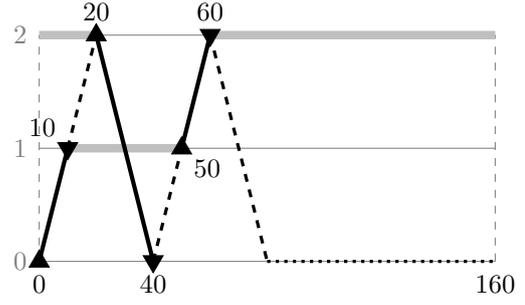
\begin{figure}
\centering
\begin{subfigure}[b]{0.45\textwidth}
\centering
\begin{tikzpicture}
\def\xmax{200}
\def\ymax{2}
\def\ctime{200}
\def\W{6cm}
\def\H{3cm}

\begin{scope}[x=\W/\xmax, y=\H/\ymax]

\foreach \h in {0,1,2}{
  \draw[dTank] (0,\h) -- (\ctime,\h);
  \node[nTank, left] at (-1,\h) {\h};
}
\draw[dTank,dashed] (0,0) -- (0,2);
\draw[dTank,dashed] (\ctime,0) -- (\ctime,2);

\coordinate (p0) at ( 0,0);
\coordinate (p1) at (50,1);
\coordinate (p2) at (180,2);
\coordinate (d0) at (200,0);
\coordinate (d1) at (10,1);
\coordinate (d2) at (60,2);
\coordinate (r0) at (80,0);

\draw[dProcess] (d1) -- (p1);
\draw[dProcess] (d2) -- (p2);

\draw[dMove] (p0) -- (d1);
\draw[dMove] (p1) -- (d2);
\draw[dMove] (p2) -- (d0);
\draw[dWait] (d1) -- (p1);
\draw[dWait] (d2) -- (p2);

\node[nPick] at (p0) {};
\node[nPick] at (p1) {};
\node[nPick] at (p2) {};
\node[nDrop]   at (d0) {};
\node[nDrop]   at (d1) {};
\node[nDrop]   at (d2) {};

\node[nTime,below=2pt of p0] {0};
\node[nTime,above=2pt of d1] {10};
\node[nTime,below=2pt of p1] {50};
\node[nTime,above=2pt of d2] {60};
\node[nTime,above=2pt of p2] {180};
\node[nTime,below=2pt] at (\ctime,0) {200};
\end{scope}
\end{tikzpicture}
\caption{Schedule ($C = 200$) obtained by the formulations of \citet{phillips1976mathematical} and \citet{leung2004optimal}.}
\label{fig:cntex1:pu}
\end{subfigure}
\hfill
\begin{subfigure}[b]{0.45\textwidth}
\centering
\begin{tikzpicture}
\def\xmax{160}
\def\ymax{2}
\def\ctime{160}
\def\W{6cm}
\def\H{3cm}

\begin{scope}[x=\W/\xmax, y=\H/\ymax]

\foreach \h in {0,1,2}{
  \draw[dTank] (0,\h) -- (\ctime,\h);
  \node[nTank, left] at (-1,\h) {\h};
}
\draw[dTank,dashed] (0,0) -- (0,2);
\draw[dTank,dashed] (\ctime,0) -- (\ctime,2);

\coordinate (p0) at ( 0,0);
\coordinate (p1) at (50,1);
\coordinate (p2) at (20,2);
\coordinate (d0) at (40,0);
\coordinate (d1) at (10,1);
\coordinate (d2) at (60,2);
\coordinate (r0) at (80,0);

%\draw[dProcess] (d0) -- (\ctime,0);
\draw[dProcess] (0,2) -- (p2);
\draw[dProcess] (d2) -- (\ctime,2);
\draw[dProcess] (d1) -- (p1);

\draw[dMove] (p0) -- (d1);
\draw[dMove] (p1) -- (d2);
\draw[dMove] (p2) -- (d0);
\draw[dTravel] (d1) -- (p2);
\draw[dTravel] (d0) -- (p1);
\draw[dTravel] (d2) -- (r0);
\draw[dWait] (r0) -- (\ctime,0);

\node[nPick] at (p0) {};
\node[nPick] at (p1) {};
\node[nPick] at (p2) {};
\node[nDrop]   at (d0) {};
\node[nDrop]   at (d1) {};
\node[nDrop]   at (d2) {};

\node[nTime,below=2pt of p0] {0};
\node[nTime,below=2pt of d0] {40};
\node[nTime,above left= 1pt and 1pt of d1] {10};
\node[nTime,below right= 1pt and 1pt of p1] {50};
\node[nTime,above=2pt of p2] {20};
\node[nTime,above=2pt of d2] {60};
\node[nTime,below=2pt] at (\ctime,0) {160};
\end{scope}
\end{tikzpicture}
\caption{
Optimal schedule ($C = 160$).
In order to meet the minimum soak time of operation~2, the hoist waits at the load station at the end of the cycle.
}
\label{fig:cntex1:opt}
\end{subfigure}
\caption{
Example \texttt{ex1}.
Formulations that are restricted with respect to the cycle finish may lead to suboptimal solutions.
}
\end{figure}

Consider the mini example with two soaking operations indicated in~\Cref{tab:cntex1}.
Note that the load and unload stations are associated, but loading and unloading constraints are unnecessary as $L_0 = L_{n+1} = 0$ and $U_0 = U_{n+1} = \infty$.
There are two possible orders of moves: $(0,1,2)$ and $(0,2,1)$.
In the first case, the minimum cycle time is $d_0 + L_1 + d_1 + L_2 + d_3 = 200$, see \Cref{fig:cntex1:pu}.
In the latter case, the cycle time is at least $d_2 + e_{0,1} + d_1 + L_2 = 160$, which can be reached, see \Cref{fig:cntex1:opt}.
The situation changes when the hoist must immediately return to the load station after executing the latest move, and then the cycle ends.
In case of move order $(0,2,1)$, move~$2$ cannot start earlier than $L_2 - e_{2,0} = 100$, but move~$1$ cannot start later than $d_0 + U_1 = 110$, which is impossible.
Thus, for the model of \citet{phillips1976mathematical}, only the first move order is feasible, so it results in a suboptimal solution.

\paragraph{Multifunction tanks.}
Note that constraints~(\ref{eq:pu:overc})-(\ref{eq:pu:overji}) do not prevent operations using the same tank from overlapping.
Recall that constraint~(\ref{eq:pu:overij}) aims to enforce that if $\lpvar{t}_{i-1} < \lpvar{t}_{j-1}$, then $\lpvar{t}_{i} < \lpvar{t}_{j-1}$.
As \citet{liu2002cyclic} pointed out, this constraint does not prevent the case $\lpvar{t}_{i-1} + d_{i-1} < \lpvar{t}_j < \lpvar{t}_i < \lpvar{t}_{j-1} + d_{j-1}$, that is, operation~$j$ crosses the cycle boundary and overlaps with operation~$i$.
The authors provided an example, in which the model produces an incorrect solution.

We also mention that constraint~(\ref{eq:pu:overij}) may also fail to prevent the case $\lpvar{t}_{i-1} < \lpvar{t}_{j-1} < \lpvar{t}_i$.
The constraint only forces a minimum gap between $\lpvar{t}_{i-1}$ and $\lpvar{t}_{j-1}$ based on the lower bound~$L_i$, however, the operation may have a duration longer than~$L_i$, and thus longer than the enforced gap.
Thus, $\lpvar{t}_i$ may be greater than $\lpvar{t}_{j-1}$.

\subsection{Formulation of \citet{leung2004optimal}}
\label{sec:shcsp:leung}
\begin{subequations}
\renewcommand{\theequation}{B\arabic{equation}}

\citet{leung2004optimal} revised the model of \citet{phillips1976mathematical} and proposed a tighter formulation, where two big-M constraints of the former model are simplified, and several valid inequalities is introduced.

We use the cycle time variable~$\lpvar{C}$, even though the authors only introduced an auxiliary variable along with the constraint $\lpvar{d} = \sum_{i=1}^n ( d_i + e_{i+1,0} ) \lpvar{z}_i$.
Again, we note that the formulation remains equivalent to the original one.

Constraints (\ref{eq:pu:cvar})-(\ref{eq:pu:2}) are used unchanged, see (\ref{eq:leung:cvar})-(\ref{eq:leung:2}).
\begin{eqnarray}
\label{eq:leung:cvar} \lpvar{C} = \lpvar{t}_{\max} + \sum_{i=1}^n ( d_i + e_{i+1,0} ) \lpvar{z}_i && \\
\label{eq:leung:3} \sum_{i=1}^n \lpvar{z}_i = 1 && \\
\label{eq:leung:1} \lpvar{t}_{\max} \geq \lpvar{t}_i && 1 \leq i \leq n\\
\label{eq:leung:2} \lpvar{t}_{\max} \leq \lpvar{t}_i + \bigM(1 - \lpvar{z}_i) && 1 \leq i \leq n
\end{eqnarray}
Travel time constraints are separated for tank~0 (see constraint~(\ref{eq:leung:4})) and for tanks $1,\ldots,n$ (see constraints~(\ref{eq:leung:5})-(\ref{eq:leung:6})). That is, $\lpvar{y}_{0,j}$ is explicitly substituted by~1.
\begin{eqnarray}
\label{eq:leung:4} \lpvar{t}_j \geq d_0 + e_{1,j} && 1 \leq j \leq n\\
\label{eq:leung:5} \lpvar{t}_j \geq (\lpvar{t}_i + d_i) + e_{i+1,j} - \bigM(1-\lpvar{y}_{i,j}) && 1 \leq i < j \leq n\\
\label{eq:leung:6} \lpvar{t}_i \geq (\lpvar{t}_j + d_j) + e_{j+1,i} - \bigM\lpvar{y}_{i,j} && 1 \leq i < j \leq n 
\end{eqnarray}
The authors observed that the big-M terms can be omitted from the right-hand-side of the soaking constraint~(\ref{eq:leung:910}) and from the left-hand-side of the soaking constraint~(\ref{eq:leung:1112}).
\begin{eqnarray}
\label{eq:leung:8} L_1 \leq \lpvar{t}_1 - d_0 \leq U_1 && \\
\label{eq:leung:910} L_i - \bigM(1-\lpvar{y}_{i-1,i}) \leq \lpvar{t}_i - (\lpvar{t}_{i-1} + d_{i-1}) \leq U_i && 2 \leq i \leq n \\
\label{eq:leung:1112} L_i \leq (\lpvar{C} + \lpvar{t}_i) - (\lpvar{t}_{i-1} + d_{i-1}) \leq U_i + \bigM\lpvar{y}_{i-1,i} && 2 \leq i \leq n
\end{eqnarray}

\citet{leung2004optimal} also proposed valid inequalities to tighten their formulation.
These constraints are based on the observation, that the latest move has no successors and all other moves precede it (see constraints~(\ref{eq:leung:16}) and~(\ref{eq:leung:17})), while a non-latest move has at least one successor (see constraint~(\ref{eq:leung:18})).
Note that the first two inequalities complement each other, that is, the inequality $1 - \lpvar{y}_{j,i} \leq 1 - \lpvar{z}_i$ ($1 \leq j < i \leq n$) is the same as constraint~(\ref{eq:leung:17}).
\begin{eqnarray}
\label{eq:leung:16} \lpvar{y}_{i,j} \leq 1 - \lpvar{z}_i && 1 \leq i < j \leq n\\
\label{eq:leung:17} \lpvar{z}_j \leq \lpvar{y}_{i,j} && 1 \leq i < j \leq n\\
\label{eq:leung:18} 1 - \lpvar{z}_i \leq \sum_{j=1}^{i-1}(1-\lpvar{y}_{j,i}) + \sum_{j=i+1}^{n}\lpvar{y}_{i,j} && 1 \leq i \leq n
\end{eqnarray}
\end{subequations}

In summary, the formulation of \citet{leung2004optimal} aims to minimize~$\lpvar{C}$ subject to constraints~(\ref{eq:leung:cvar})-(\ref{eq:leung:18}), where
$\lpvar{t}_0 = 0$, $0\leq \lpvar{t}_i$ ($1\leq i\leq n$), $\lpvar{y}_{i,j} \in \{0,1\}$ ($1\leq i < j \leq n$), $\lpvar{z}_i \in \{0,1\}$ ($1\leq i \leq n$).

\subsubsection{Remarks, issues}
In the following, we discuss our observations regarding the formulation.

\paragraph{Unspecified loading and unloading constraints.}
Similarly to \citet{phillips1976mathematical}, \citet{leung2004optimal} also did not specify the loading and unloading constraints for the associated load and unload stations.

\paragraph{Narrowed solution space.}
Similarly to \citet{phillips1976mathematical}, the formulation of \citet{leung2004optimal} is also restricted with respect to the cycle finish.

\subsection{Formulation of \citet{zhou2003single}}\label{sec:shcsp:zhou}

Constraints (\ref{eq:pu:3})-(\ref{eq:pu:89}) are used unchanged, see constraints~(\ref{eq:zhou:2})-(\ref{eq:zhou:711}).
The cycle time variable~$\lpvar{C}$ is now a lower bound on the return time, see constraint~(\ref{eq:zhou:1}).
Loading and unloading are now explicitly taken into account; constraint~(\ref{eq:zhou:910}) refers to associated load and unload configuration (cf. constraint~(\ref{eq:sc:load:ass})).

\begin{subequations}
\renewcommand{\theequation}{C\arabic{equation}}
\begin{eqnarray}
\label{eq:zhou:1} \lpvar{C} \geq \lpvar{t}_{\max} + \sum_{i=1}^n ( d_i + e_{i+1,0} ) \lpvar{z}_i && \\
\label{eq:zhou:2} \sum_{i=1}^n \lpvar{z}_i = 1 && \\
\label{eq:zhou:3} \lpvar{t}_{\max} \geq \lpvar{t}_i && 1 \leq i \leq n\\
\label{eq:zhou:4} \lpvar{t}_{\max} \leq \lpvar{t}_i + \bigM(1 - \lpvar{z}_i) && 1 \leq i \leq n\\
\label{eq:shscp:zhou:5}\lpvar{t}_j \geq (\lpvar{t}_i + d_i) + e_{i+1,j} - \bigM(1-\lpvar{y}_{i,j}) && 0 \leq i < j \leq n\\
\label{eq:shscp:zhou:6}\lpvar{t}_i \geq (\lpvar{t}_j + d_j) + e_{j+1,i} - \bigM\lpvar{y}_{i,j} && 1 \leq i < j \leq n:\ i + 1 < j\\
\label{eq:zhou:812}L_i - \bigM(1-\lpvar{y}_{i-1,i}) \leq \lpvar{t}_i - (\lpvar{t}_{i-1} + d_{i-1}) \leq U_i + \bigM(1-\lpvar{y}_{i-1,i}) && 1 \leq i \leq n\\
\label{eq:zhou:711}
L_i - \bigM\lpvar{y}_{i-1,i} \leq (\lpvar{C} + \lpvar{t}_i) - (\lpvar{t}_{i-1} + d_{i-1}) \leq U_i + \bigM\lpvar{y}_{i-1,i} && 1 \leq i \leq n\\
\label{eq:zhou:910} L_0 \leq \lpvar{C} - (\lpvar{t}_n + d_n) \leq U_0&&
\end{eqnarray}

In summary, the formulation of \citet{zhou2003single} aims to minimize~$\lpvar{C}$ subject to constraints~(\ref{eq:zhou:1})-(\ref{eq:zhou:910}), where
$\lpvar{t}_0 = 0$, $0\leq \lpvar{t}_i$ ($1\leq i\leq n$), $\lpvar{y}_{i,j} \in \{0,1\}$ ($1\leq i < j \leq n$).

\subsubsection{Extension to multi-tank operations}
\citet{zhou2003single} extended their formulation to the multi-tank case.
For a bottleneck operation~$i$, the corresponding constraints~(\ref{eq:zhou:812})-(\ref{eq:zhou:711}) are replaced with constraints (\ref{eq:zhou:14})-(\ref{eq:zhou:13}), (cf. constraint (\ref{eq:msoak:M})-(\ref{eq:msoakw:M})).
\begin{eqnarray}
\label{eq:zhou:14}L_i - \bigM(1-\lpvar{y}_{i-1,i}) \leq ((m-1)\lpvar{C} + \lpvar{t}_i) - (\lpvar{t}_{i-1} + d_{i-1}) \leq U_i + \bigM(1-\lpvar{y}_{i-1,i}) &&\\
\label{eq:zhou:15}
L_i - \bigM\lpvar{y}_{i-1,i} \leq (m\lpvar{C} + \lpvar{t}_i) - (\lpvar{t}_{i-1} + d_{i-1}) \leq U_i + \bigM\lpvar{y}_{i-1,i} &&\\
\label{eq:zhou:13} U_i/m \leq \lpvar{C} \leq L_i/(m-1)&&
\end{eqnarray}
\end{subequations}

\subsubsection{Remarks}
The formulation of \citet{zhou2003single} is \emph{unrestricted} with respect to the cycle finish, as the authors corrected the calculation of the cycle time with constraint~(\ref{eq:zhou:1}).
The authors also specify the loading and unloading constraints for associated load and unload stations.

\paragraph{Unspecified multifunction tank constraints.}
The authors proposed an example with a multifunction tank, however, they did not specify the corresponding constraints.

\paragraph{Invalid lower and upper bounds.}
The authors derived constraint~(\ref{eq:zhou:13}) for the multi-tank case by combining the following lower and upper bounds on the cycle time:
$L_i/m \leq \lpvar{C} \leq L_i/(m-1)$ and $U_i/m \leq \lpvar{C} \leq U_i/(m-1)$.
However, the bounds $U_i/m \leq \lpvar{C}$ and $\lpvar{C} \leq L_i/(m-1)$ are not valid in general.
Notice that a feasible solution remains feasible if the operation upper bounds are increased. 
However, in the presence of constraint~(\ref{eq:zhou:13}), an otherwise feasible problem may become infeasible if $U_i$ is increased beyond $L_i m/(m-1)$.

\subsection{Improved extended formulation}\label{sec:shcsp:imp1}
\begin{subequations}
\renewcommand{\theequation}{D\arabic{equation}}
Now, we present an improved extended formulation, which is a combination of the previous models.
Similarly to \citet{zhou2003single}, we correct the calculation of the cycle time (see constraint~(\ref{eq:imp1:cvar})).
Similarly to \citet{leung2004optimal}, we use the strengthened soaking constraints (\ref{eq:imp1:soak})-(\ref{eq:imp1:soakw}) and the valid inequalities (\ref{eq:imp1:valid1})-(\ref{eq:imp1:valid2}).

The objective is to minimize the cycle time~$\lpvar{C}$ subject to constraints (\ref{eq:imp1:cvar})-(\ref{eq:imp1:valid2}), where $0\leq \lpvar{t}_i$ ($1\leq i\leq n)$, $\lpvar{y}_{i,j} \in \{0,1\}$ ($1 \leq i,j \leq n: i\neq j$), $\lpvar{z}_i \in \{0,1\}$ ($1 \leq i \leq n$).

\begin{eqnarray}
\label{eq:imp1:cvar}\lpvar{C} \geq \lpvar{t}_{\max} + \sum_{i=1}^n \left( d_i + e_{i+1,0} \right) \lpvar{z}_i\\
\sum_{i=1}^n \lpvar{z}_i = 1 && \\
\lpvar{t}_{\max} \geq \lpvar{t}_i && 1 \leq i \leq n\\
\lpvar{t}_{\max} \leq \lpvar{t}_i + \bigM(1 - \lpvar{z}_i) && 1 \leq i \leq n\\
\lpvar{y}_{i,j} + \lpvar{y}_{j,i} = 1 && 0 \leq i < j \leq n\\
\lpvar{t}_j \geq \lpvar{t}_i + d_i + e_{i+1,j} - \bigM(1-\lpvar{y}_{i,j}) && 0 \leq i,j \leq n:\ i\neq j\\
\label{eq:imp1:soak}L_i - \bigM(1-\lpvar{y}_{i-1,i}) \leq \lpvar{t}_i - (\lpvar{t}_{i-1} + d_{i-1}) \leq U_i && 1 \leq i \leq n\\
\label{eq:imp1:soakw} L_i \leq (\lpvar{C} + \lpvar{t}_i) - (\lpvar{t}_{i-1} + d_{i-1}) \leq U_i + \bigM\lpvar{y}_{i-1,i} && 1 \leq i \leq n\\
\lpvar{t}_0 = 0 && \\
\label{eq:imp1:last}\lpvar{y}_{0,j} = 1 && 1 \leq j \leq n
\end{eqnarray}
The following are valid inequalities:
\begin{eqnarray}
\label{eq:imp1:valid1} \lpvar{y}_{i,j} \leq 1 - \lpvar{z}_i && 1 \leq i, j \leq n: i \neq j \\
\label{eq:imp1:valid2} 1 - \lpvar{z}_i \leq \sum_{\substack{j=1\\ j\neq i}}^{n}\lpvar{y}_{i,j} && 1 \leq i \leq n
\end{eqnarray}
\end{subequations}
Depending on the configuration of the load and unload stations, constraint~(\ref{eq:sc:load:diss}) or (\ref{eq:sc:load:ass}) can be used.
In case of multifunction tanks, constraint~(\ref{eq:sc:multi}); in case of multi-tank operations, constraints~(\ref{eq:mstage:sum}) and (\ref{eq:mstage:soak})-(\ref{eq:mstage:soakw}) are added to the model.

\subsection{Formulation of \citet{liu2002cyclic}}\label{sec:shcsp:liu}
\citet{liu2002cyclic} considered the scenario in which the values~$d_i$ are only lower bounds on the move times, rather than exact values.
For this, the non-negative continuous \emph{slack time variable}~$\lpvar{d}^+_i$ measures the slack time between the lower bound~$d_i$ and the realized time for move~$i$ ($0 \leq i \leq n$).
By this, the start time of operation~$i$ ($1\leq i\leq n$) is $\lpvar{t}_{i-1} + d_{i-1} + \lpvar{d}^+_{i-1}$.
Note that fixing $\lpvar{d}^+_i = 0$ for each $0\leq i\leq n$ results in a base formulation with exact move times.

\begin{subequations}
\renewcommand{\theequation}{E\arabic{equation}}
The objective is to minimize the cycle time~$\lpvar{C}$, where the cycle time is long enough for the hoist to return to the load station after executing any move, see constraint~(\ref{eq:liu:8}).
\begin{eqnarray}
\label{eq:liu:8} \lpvar{C} \geq \lpvar{t}_i + d_i + \lpvar{d}^+_i + e_{i+1,0} && 1 \leq i \leq n
\end{eqnarray}
Constraint~(\ref{eq:liu:2}) is used only if the load and unload stations are associated.
Constraints~(\ref{eq:liu:3})-(\ref{eq:liu:4}) ensure that the soak time bounds are respected.
\begin{eqnarray}
L_0 \leq \label{eq:liu:2}\lpvar{C} - (\lpvar{t}_n + d_n + \lpvar{d}^+_n) \leq U_0 && \text{if }s_{n+1} = s_0\\
\label{eq:liu:3}L_i - \bigM(1-\lpvar{y}_{i-1,i}) \leq \lpvar{t}_i - (\lpvar{t}_{i-1} + d_{i-1} + \lpvar{d}^+_{i-1}) \leq U_i + \bigM(1-\lpvar{y}_{i-1,i}) && 1 \leq i \leq n\\
\label{eq:liu:4}L_i - \bigM\lpvar{y}_{i-1,i} \leq (\lpvar{C} + \lpvar{t}_i) - (\lpvar{t}_{i-1} + d_{i-1} + \lpvar{d}^+_{i-1}) \leq U_i + \bigM\lpvar{y}_{i-1,i} && 1 \leq i \leq n
\end{eqnarray}
Constraints~(\ref{eq:liu:5})-(\ref{eq:liu:7}) refer to the travel time constraints.
\begin{eqnarray}
\label{eq:liu:5}\lpvar{t}_j \geq d_0 + \lpvar{d}^+_0 + e_{1,j} && 1 \leq j \leq n\\
\label{eq:liu:6}\lpvar{t}_j \geq \lpvar{t}_i + d_i + \lpvar{d}^+_i + e_{i+1,j} - \bigM(1-\lpvar{y}_{i,j}) && 1 \leq i < j \leq n\\
\label{eq:liu:7}\lpvar{t}_i \geq \lpvar{t}_j + d_j + \lpvar{d}^+_j + e_{j+1,i} - \bigM\lpvar{y}_{i,j} && 1 \leq i < j \leq n
\end{eqnarray}

In summary, the formulation of~\citet{liu2002cyclic} aims to minimize~$\lpvar{C}$ subject to constraints~(\ref{eq:liu:8})-(\ref{eq:liu:7}), where $0\leq \lpvar{t}_i$ ($1\leq i\leq n$), $\lpvar{y}_{i,j} \in \{0,1\}$ ($1\leq i < j \leq n$), $0\leq \lpvar{d}^+_i$ ($1\leq i\leq n$).

\subsubsection{Extension to multifunction tanks}
\citet{liu2002cyclic} extended their model to multifunction tanks with the following constraints (cf. \Cref{sec:model:multi}).
\begin{eqnarray}
\label{eq:liu:12}\lpvar{y}_{i-1,i} + \lpvar{y}_{j-1,j} \geq 1 && 0 \leq i + 1 < j \leq n:\ s_i = s_j\\
\label{eq:liu:13}\lpvar{y}_{i,j-1} + (1 - \lpvar{y}_{i-1,j}) \geq 3 - (\lpvar{y}_{i-1,i} + \lpvar{y}_{j-1,j}) && 0 \leq i + 1 < j \leq n:\ s_i = s_j
\end{eqnarray}
\end{subequations}

\subsubsection{Remarks}
\citet{liu2002cyclic} presented a base formulation by neglecting additional variables~$\lpvar{t}_{\max}$ and $\lpvar{z}_i$.
Consequently, the formulation is unrestricted with respect to the cycle finish, see constraint~(\ref{eq:liu:8}).
The authors explicitly specified that the loading constraints are used only for associated load and unload stations, see constraint~(\ref{eq:liu:2}).
The authors corrected the multifunction tank constraints of~\citet{phillips1976mathematical}.
However, note that constraint~(\ref{eq:liu:12}) is redundant in the presence of constraint~(\ref{eq:liu:13}), cf. constraint~(\ref{eq:sc:multi}).

\subsection{Improved base formulation}\label{sec:shcsp:imp2}
\begin{subequations}
\renewcommand{\theequation}{F\arabic{equation}}
Now, we present a strengthened version of the base formulation of \citet{liu2002cyclic}, where the stronger soaking constraints are used.
The objective is to minimize the cycle time~$\lpvar{C}$ subject to constraints (\ref{eq:imp2:cvar})-(\ref{eq:imp2:last}), where $0\leq \lpvar{t}_i$ ($1\leq i\leq n)$, $\lpvar{y}_{i,j} \in \{0,1\}$ ($1 \leq i,j \leq n: i\neq j$).

\begin{eqnarray}
\label{eq:imp2:cvar}\lpvar{C} \geq \lpvar{t}_i + d_i + e_{i+1,0} && 1 \leq i \leq n\\
\lpvar{t}_j \geq \lpvar{t}_i + d_i + e_{i+1,j} - \bigM(1-\lpvar{y}_{i,j}) && 1 \leq i, j \leq n:\ i\neq j\\
\lpvar{y}_{i,j} + \lpvar{y}_{j,i} = 1 && 1 \leq i < j \leq n\\
L_i - \bigM(1-\lpvar{y}_{i-1,i}) \leq \lpvar{t}_i - (\lpvar{t}_{i-1} + d_{i-1}) \leq U_i && 1 \leq i \leq n\\
L_i \leq (\lpvar{C} + \lpvar{t}_i) - (\lpvar{t}_{i-1} + d_{i-1}) \leq U_i + \bigM\lpvar{y}_{i-1,i} && 1 \leq i \leq n\\
\lpvar{t}_0 = 0 && \\
\label{eq:imp2:last}\lpvar{y}_{0,j} = 1 && 1 \leq j \leq n
\end{eqnarray}
\end{subequations}
Depending on the configuration of the load and unload stations, constraint~(\ref{eq:sc:load:diss}) or (\ref{eq:sc:load:ass}) can be used.
In case of multifunction tanks, constraint~(\ref{eq:sc:multi}); in case of multi-tank operations, constraints~(\ref{eq:mstage:sum}) and (\ref{eq:mstage:soak})-(\ref{eq:mstage:soakw}) are added to the model.

\section{Computational experiments}\label{sec:comp}
In this section, we present the results of our computational experiments to compare the proposed formulations. 
We evaluate them in three complementary ways: 
by examining the strength of the MIP~formulations via their LP-relaxations (\Cref{sec:comp:relax}); 
by testing their performance in solving the MIP models on a single thread (\Cref{sec:comp:mip1}); 
and by assessing their multi-thread performance on 8~threads (\Cref{sec:comp:mip8}). 
These three perspectives allow us to isolate the effect of the formulation itself (LP-relax), the solver heuristics and branching (single-thread MIP), and parallelization (multi-thread MIP) on solution quality and computational efficiency.

Note that the proposed CP~formulation solves all considered instances to optimality within short running times, and thus provides a simple and effective baseline for practical use.
In contrast, MIP formulations enable a deeper analysis of modeling choices, such as the impact of constraint structure and formulation strength on LP relaxations and solver performance.
Therefore, the two approaches play complementary roles in our study.

In the following tables, we report only summary results.
%In the following tables, we report only summary results, and refer the reader to the supplementary material for the detailed results.
Column \textit{Phillips} refers to the formulation (\ref{eq:pu:cvar})-(\ref{eq:pu:89}) of \citet{phillips1976mathematical}.
Column \textit{Leung} refers to the formulation (\ref{eq:leung:cvar})-(\ref{eq:leung:1112}) of \citet{leung2004optimal}.
In case of \textit{Leung$^+$}, the model is also extended with the valid inequalities (\ref{eq:leung:16})-(\ref{eq:leung:18}).
Column \textit{Zhou} refers to the formulation (\ref{eq:zhou:1})-(\ref{eq:zhou:910}) of \citet{zhou2003single}.
Column \textit{Imp1} refers to the improved extended formulation (\ref{eq:imp1:cvar})-(\ref{eq:imp1:last}) proposed in \Cref{sec:shcsp:imp1}.
In case of \textit{Imp$^+$}, the model is also extended with the valid inequalities (\ref{eq:imp1:valid1})-(\ref{eq:imp1:valid2}).
Column \textit{Liu} refers to the formulation (\ref{eq:liu:2})-(\ref{eq:liu:13}) of \citet{liu2002cyclic} with exact move times (i.e., $\lpvar{d}^+_i = 0$ is fixed for $0\leq i\leq n$).
Column \textit{Imp2} refers to the improved base formulation (\ref{eq:imp2:cvar})-(\ref{eq:imp2:last}) proposed in \Cref{sec:shcsp:imp2}.
In each case, the objective is to minimize the cycle time~$\lpvar{C}$.

In the comparison, we take into account that the formulations differ in their main properties, summarized in~\Cref{tab:shcsp:form}.
Note that the loading and unloading constraints were specified differently (or not specified at all) in the formulations, but all formulations could handle both constraints~(\ref{eq:sc:load:diss}) and~(\ref{eq:sc:load:ass}).
In order to make a fair comparison, we generated instances with~$L_0 = U_0 = 0$ and $L_{n+1} = U_{n+1} = \infty$,
thus loading and unloading constraints were irrelevant.

\begin{table}
\centering
\begin{threeparttable}[b]
\caption{Main properties of the formulations.}
\label{tab:shcsp:form}
\begin{tabular}{lcccc}
\toprule
Formulation & Type & Cycle finish & Soaking constraints & Valid inequalities \\
\midrule
Phillips  & extended & restricted   & original     & no \\
Leung     & extended & restricted   & strengthened & no \\
Leung$^+$ & extended & restricted   & strengthened & yes\\
Zhou      & extended & unrestricted & original     & no \\
Imp1      & extended & unrestricted & strengthened & no \\
Imp1$^+$  & extended & unrestricted & strengthened & yes\\
Liu       & base     & unrestricted & original     & no \\
Imp2      & base     & unrestricted & strengthened & no \\
\bottomrule
\end{tabular}
\end{threeparttable}
\end{table}

\paragraph{Setup.}
The models are implemented in Python (version 3.13.11) using Google OR-Tools (version 9.14.6206).
The CP models are solved with OR-Tools CP-SAT, while the MIP models are solved with Gurobi (version 12.0.3).
All experiments were conducted in a cloud environment on an Ubuntu 22.04 system equipped with 16 single-threaded Intel Xeon Cascade Lake cores running at 2.1 GHz.

\subsection{Instances, preliminary experiments}\label{sec:comp:inst}

In \Cref{sec:comp:inst:bench}, we present the widely-used benchmark instances that can be found in the literature.
There, we also explain the inconsistencies between the optimal values reported, see~\Cref{tab:indicated}.
In~\Cref{sec:comp:inst:rand}, we present the procedure for generating random instances that can challenge the modern CP and MIP~solvers.

\subsubsection{Benchmark instances}\label{sec:comp:inst:bench}

\citet{phillips1976mathematical} proposed an example (PU) for which they reported a solution with a cycle time of~580.
Contrary to their claim, this hoist schedule is not optimal; the true optimal cycle time is~521.
The "optimal" cycle time reported by \citet{liu2025fast} suggests that either their implementation differs from the intended formulation, or a modified instance was used.
Our instance \texttt{philu} corresponds to the one proposed by \citet{phillips1976mathematical}, converted into an equivalent instance with dissociative load and unload stations (see \Cref{sec:shcsp:pre}).
This modification applies to all the following instances as well.

According to \citet{leung2004optimal, laajili2021adapted}, instance Mini is the truncation of instance PU to the first eight tanks.
However, the instance itself is not provided, which may cause inconsistencies.
For example, it is not clear whether further adjustment of move times is needed, since \citet{phillips1976mathematical} defined their original instance with $d_i = e_{i,i+1} + 20$.
Our instance \texttt{philu\_mini} is made without any further adjustment.
According to our CP model, the optimal cycle time is~284 for dissociated load and unload stations (i.e., with constraint~(\ref{eq:sc:load:diss})), and the optimum is~340 for associated load and unload stations (i.e., with constraint~(\ref{eq:sc:load:ass})).

\citet{shapiro1988hoist} applied their solution approach on four new test problems, Black Oxide~1 (BO1), Black Oxide~2 (BO2), Zinc (Zn), and Copper (Cu), which were provided by a vendor of wetlines.
The data for these examples were not provided in their paper, but later \citet{leung2004optimal} proposed these or similar instances.
The original instances included duplicated tanks, allowing for multi-tank operations, however, the later instances are typically solved only with single-tank operations.
This modeling difference leads to discrepancies in the reported cycle times, as the underlying problem structures are not fully equivalent.
Our instances \texttt{bo1}, \texttt{bo2}, \texttt{zn}, and \texttt{cu} are based on \citep{laajili2021adapted}, but the time values have been multiplied by 10 to obtain integer values, as the CP-SAT solver operates on integers.
We emphasize that this transformation does not change the set of feasible solutions.
It can be verified, for example, using our CP model, that the optimal objective value for BO1 changes depending on the configuration of the load and unload stations (dissociated: 281.9; associated: 333.2).
Note that we could not reproduce the "optimal" cycle time~299.5 reported in \citep{leung2004optimal}.
For the other three instances, there is no difference in these two cases (BO2: 279.3, Cu: 1847.2, Zn: 1743.4).
\citet{liu2025fast} probably rounded the time values to obtain integer values for CP-SAT, however, we could not reproduce their reported cycle times for these instances.

\citet{manier1994contribution} proposed two new industrial instances (Ligne~1 and Ligne~2), also available in \citep{manier2008evolutionary}.
Ligne~1 originally included a circulation constraint, therefore there was a separate unloading operation associated with a non-zero lower bound.
\citep{laajili2021adapted} slightly modified this instance, removing the unloading operation and adding its lower bound to the lower bound of the loading operation.
Our instances \texttt{ligne1} and \texttt{ligne2} are the ones provided in \citep{laajili2021adapted}.
It can be verified, for example, using our CP model, that the optimal cycle time for Ligne~1 depends on the load-unload configuration (dissociated: 392; associated: 425).
In case of Ligne~2, the optimal cycle time is 712 in both cases.

\subsubsection{Random instances}\label{sec:comp:inst:rand}
Some papers proposed procedures for generating random instances, but without providing the instances themselves.
Our preliminary experiments showed that the method of \citet{che2015robust} generates instances that can challenge the modern CP and MIP solvers, and therefore suitable for our comparisons.
The authors generated instances for a robust, single-hoist cyclic scheduling problem with dissociated load and unload stations with parameters~$n \in \{12,14,16,18,20,22\}$ and~$\mu \in \{1.5, 2.0, 2.5\}$ in the following way: $e_{i,i+1} = 1 + \operatorname{uniform}(0,4)$ ($0\leq i\leq n$); $e_{i,j+1} = e_{j+1,i} = \sum_{k=i}^j e_{k,k+1}$ ($1\leq i < j \leq n$); $d_i = e_{i,i+1}+12$, $L_i = 40 + \operatorname{uniform}(0,140)$ and $U_i = [\mu L_i]$ ($0\leq i\leq n$), where $\operatorname{uniform}(a,b)$ refers to a uniformly random integer from $\{a,a+1,\ldots,b\}$.

Our preliminary experiments showed that, some of the instances generated in this way are infeasible in the presence of loading and unloading constraints, especially for the formulations of \citet{phillips1976mathematical, leung2004optimal} that are restricted with respect to the cycle finish.
Thus, we slightly modified this procedure by setting $L_0 = 0$ and $U_0 = \infty$.
The implemented generation procedure can be found in our library.
For our experiments, we generated 10~instances for each pair of parameters from $n \in \{14,19,24\}$ and $\mu \in \{1.5,2.0,2.5\}$, resulting in 90 instances.
Note that for a fixed~$n$, the $k$th instances for the different $\mu$ parameters differ only in the $U_i$ values, that is, the empty travel time matrix, the move times, and the lower bounds are the same.

\subsection{LP-relaxations}\label{sec:comp:relax}

We computed the root LP bound for each MIP~formulation by relaxing all integer and binary variables to continuous within their original bounds (0–1 for binaries), and disabling presolve (\texttt{presolve = Emphasis.OFF}).
\Cref{tab:shcsp:lp1,tab:shcsp:lp2} indicate the computational results on the benchmark instances, and summarize the average computational results on the randomly generated instances, respectively.
Gurobi solved each instance in a split of seconds, thus solution times are not indicated.
Recall that comparing the restricted formulations of \citet{phillips1976mathematical, leung2004optimal} to the others may not be fair, since the former formulations have a narrower feasible solution space.

\paragraph{A note on the operation upper bounds.}
Recall that in our instance generation procedure, the upper bounds on the soaking times were scaled as $U_i=\mu L_i$.
Interestingly, we observed that the optimal value of the LP-relaxation remained unchanged for all tested values of $\mu$, despite the fact that the feasible region of the integer problem was modified.
Thus, for these types of instances, the upper bounds~$U_i$ had little effect in the relaxation.

\paragraph{Strengthened constraints, valid inequalities.}
The significant improvement of \textit{Leung} over \textit{Phillips} came from the use of stronger soaking constraints.
Using additional valid inequalities on top of this had only a tiny effect on LP-relaxation (compare \textit{Leung} and \textit{Leung$^+$}).
The same is observed for the unrestricted formulations.
Using stronger soaking constraints strengthened the formulations (compare \textit{Zhou} and \textit{Imp1}; \textit{Liu} and \textit{Imp2}), however, the valid inequalities did not provide much benefit when strong soaking constraints were already in place (compare \textit{Imp1} and \textit{Imp1$^+$}).
In case of some industrial instances (\texttt{cu}, \texttt{zn}), the difference is significant between base and stronger soaking constraints.

\paragraph{Base formulations vs. Extended formulations.}
The base formulations always provided better or equal LP-relaxations than the corresponding extended formulations.
Considering the original soaking constraints, the base formulation (\textit{Liu}) provided 13.9\% stronger LP bounds on average than the extended formulations (\textit{Phillips}, \textit{Zhou}).
The average improvement was 6.5\%, 9.3\% and 27.0\% on the instances with $n=14$, $19$, and $24$, respectively.
The improvement on a single instance was up to $63.8\%$.
In the case of strengthened soaking constraints, the base formulation (\textit{Imp2}) provided 7.1\% better LP bounds on average than the extended formulations (\textit{Leung}, \textit{Leung$^+$}, \textit{Imp1}, \textit{Imp1$^+$}).
The average improvement was 3.2\%, 7.3\% and 11.0\% on the instances with $n=14$, $19$, and $24$, respectively.
The improvement on a single instance was up to $24.5\%$.

\begin{table}
\centering
\begin{threeparttable}[b]
\caption{LP-relaxation values on the benchmark instances.}
\label{tab:shcsp:lp1}
\begin{tabular}{l r r rr r rr r r}
\toprule
 && \multicolumn{8}{c}{Formulations}\\
\cmidrule(lr){3-10}
\multicolumn{1}{c}{Instance} & OPT\tnote{*} & Phillips & Leung & Leung$^+$ & Zhou & Imp1 & Imp1$^+$ & Liu & Imp2\\
\midrule
\texttt{philu}  &   521 & 211.0  &  211.0 &  211.0 &  211.0 &  211.0 &  211.0  &  217.0 &  217.0 \\
\texttt{philu\_mini}  &   284 & 203.0 & 203.0 & 203.0 & 203.0 & 203.0 & 203.0 & 217.0 & 217.0 \\
\texttt{cu}     & 18472 & 2756.4 & 9316.5 & 9316.6 & 2756.0 & 9316.6 & 9316.6  & 2834.0 & 9404.5 \\
\texttt{zn}     & 17434 & 2709.0 & 8725.5 & 8725.6 & 2709.0 & 8725.5 & 8725.5  & 2875.0 & 8879.0 \\
\texttt{bo1}    &  2819 & 2205.0 & 2205.0 & 2206.3 & 2205.5 & 2205.0 & 2206.3  & 2378.0 & 2378.0 \\
\texttt{bo2}    &  2793 & 2205.0 & 2205.0 & 2206.5 & 2205.0 & 2205.0 & 2206.5  & 2378.0 & 2378.0 \\
\texttt{ligne1} &   392 &  107.0 &  188.0 &  188.0 &  107.0 &  188.0 &  188.0  &  124.0 &  191.0 \\
\texttt{ligne2} &   712 &  215.0 &  369.0 &  369.1 &  215.0 &  369.0 &  369.1  &  224.0 &  380.5 \\
\bottomrule
\end{tabular}
\begin{tablenotes}
\item[*] obtained with the CP~formulation proposed in~\Cref{sec:sc:cp} without loading/unloading constraints
\end{tablenotes}
\end{threeparttable}
\end{table}

\begin{table}
\centering
\begin{threeparttable}[b]
\caption{Average LP-relaxation values on our randomly generated instances.}
\label{tab:shcsp:lp2}
\begin{tabular}{l r rr r rr r r}
\toprule
 & \multicolumn{8}{c}{Formulations}\\
\cmidrule(lr){2-9}
\multicolumn{1}{c}{Instance parameters} & Phillips & Leung & Leung$^+$ & Zhou & Imp1 & Imp1$^+$ & Liu & Imp2\\
\midrule
$n=14$, $\mu=1.5$ & 147.3 & 165.0 & 165.1 & 147.3 & 165.0 & 165.1 & 156.9 & 170.3\\
$n=14$, $\mu=2.0$ & 147.3 & 165.0 & 165.1 & 147.3 & 165.0 & 165.1 & 156.9 & 170.3\\
$n=14$, $\mu=2.5$ & 147.3 & 165.0 & 165.1 & 147.3 & 165.0 & 165.1 & 156.9 & 170.3\\
\cmidrule(lr){1-9}
$n=19$, $\mu=1.5$ & 147.9 & 156.7 & 156.7 & 147.9 & 156.7 & 156.7 & 161.6 & 168.1\\
$n=19$, $\mu=2.0$ & 147.9 & 156.7 & 156.7 & 147.9 & 156.7 & 156.7 & 161.6 & 168.1\\
$n=19$, $\mu=2.5$ & 147.9 & 156.7 & 156.7 & 147.9 & 156.7 & 156.7 & 161.6 & 168.1\\
\cmidrule(lr){1-9}
$n=24$, $\mu=1.5$ & 135.1 & 157.8 & 157.9 & 135.1 & 157.8 & 157.9 & 171.6 & 175.2\\
$n=24$, $\mu=2.0$ & 135.1 & 157.8 & 157.9 & 135.1 & 157.8 & 157.9 & 171.6 & 175.2\\
$n=24$, $\mu=2.5$ & 135.1 & 157.8 & 157.9 & 135.1 & 157.8 & 157.9 & 171.6 & 175.2\\
\midrule
$n=14$           & 147.3 & 165.0 & 165.1 & 147.3 & 165.0 & 165.1 & 156.9 & 170.3\\
$n=19$           & 147.9 & 156.7 & 156.7 & 147.9 & 156.7 & 156.7 & 161.6 & 168.1\\
$n=24$           & 135.1 & 157.8 & 157.9 & 135.1 & 157.8 & 157.9 & 171.6 & 175.2\\
\midrule
all              & 143.4 & 159.8 & 159.9 & 143.4 & 159.8 & 159.9 & 163.4 & 171.2\\
\bottomrule
\end{tabular}
\end{threeparttable}
\end{table}

\subsection{Performance tests}
For these experiments, we set a time limit of 10~minutes.
Accordingly, the status of an optimization was classified as: \textit{No solution found}, meaning that no solution was found within the time limit;
\textit{Feasible}, meaning that a feasible solution was found, but it was either not optimal or its optimality was not proved within the time limit; or
\textit{Optimal}, meaning that the optimality of the obtained solution was proved.

In the result tables, we summarize the average computational results on instances with fixed $n$ or $\mu$ parameters, and for all instances.
For these sets, we indicate
the number of instances with solution status Optimal / Feasible / No solution found;
the geometric mean of running times in seconds;
the average optimality gap, calculated as $100\% \times (\operatorname{UB}^{\operatorname{MIP}}-\operatorname{LB}^{\operatorname{MIP}})/\operatorname{UB}^{\operatorname{MIP}}$, where $\operatorname{UB}^{\operatorname{MIP}}$ is the objective value of the best found solution, and $\operatorname{LB}^{\operatorname{MIP}}$ is the best lower bound;
and the average gap to the optimum, calculated as $100\% \times (\operatorname{UB}^{\operatorname{MIP}}-\operatorname{OPT})/\operatorname{UB}^{\operatorname{MIP}}$, where $\operatorname{OPT}$ is the optimal solution.

\subsubsection{Single thread tests}\label{sec:comp:mip1}
\Cref{tab:shcsp:gurobi1} summarizes the results obtained with the Gurobi solver using a single thread (\texttt{Threads = 1}).

\begin{table}
\centering
\begin{threeparttable}[b]
\caption{Average computational results on randomly generated instances (Gurobi, single thread).}
\label{tab:shcsp:gurobi1}
\begin{tabular}{l rrrrrrrr}
\toprule
Instances
 & \multicolumn{1}{c}{Phillips}
 & \multicolumn{1}{c}{Leung}
 & \multicolumn{1}{c}{Leung$^+$}
 & \multicolumn{1}{c}{Zhou}
 & \multicolumn{1}{c}{Imp1}
 & \multicolumn{1}{c}{Imp1$^+$}
 & \multicolumn{1}{c}{Liu}
 & \multicolumn{1}{c}{Imp2}\\
\midrule
$\mu=1.5$ & 30/0/0 & 30/0/0 & 30/0/0 & 30/0/0 & 30/0/0 & 30/0/0 & 30/0/0 & 30/0/0 \\
          & 0.6 & 0.9 & 0.7 & 0.7 & 0.7 & 0.7 & 0.3 & 0.3 \\
          & 0.00\% & 0.00\% & 0.00\% & 0.00\% & 0.00\% & 0.00\% & 0.00\% & 0.00\% \\
          & 0.75\% & 0.75\% & 0.75\% & 0.00\% & 0.00\% & 0.00\% & 0.00\% & 0.00\% \\
\addlinespace[1ex]
$\mu=2.0$ & 30/0/0 & 30/0/0 & 30/0/0 & 30/0/0 & 30/0/0 & 30/0/0 & 30/0/0 & 30/0/0 \\
          & 8.8 & 8.8 & 7.3 & 8.3 & 8.5 & 6.6 & 3.7 & 3.9 \\
          & 0.00\% & 0.00\% & 0.00\% & 0.00\% & 0.00\% & 0.00\% & 0.00\% & 0.00\% \\
          & 0.00\% & 0.00\% & 0.00\% & 0.00\% & 0.00\% & 0.00\% & 0.00\% & 0.00\% \\
\addlinespace[1ex]
$\mu=2.5$ & 24/6/0 & 21/9/0 & 24/6/0 & 23/7/0 & 22/8/0 & 24/6/0 & 26/4/0 & 26/4/0 \\
          & 89.5 & 90.9 & 70.6 & 94.1 & 90.7 & 65.2 & 40.9 & 37.9 \\
          & 10.47\% & 12.43\% & 8.45\% & 11.68\% & 12.04\% & 8.37\% & 6.32\% & 6.19\% \\
          & 1.14\% & 0.80\% & 1.00\% & 1.18\% & 0.87\% & 0.67\% & 0.41\% & 0.21\% \\
\midrule
$n=14$ & 30/0/0 & 30/0/0 & 30/0/0 & 30/0/0 & 30/0/0 & 30/0/0 & 30/0/0 & 30/0/0 \\
       & 4.3 & 4.7 & 3.7 & 4.1 & 3.5 & 3.0 & 2.1 & 2.2 \\
       & 0.00\% & 0.00\% & 0.00\% & 0.00\% & 0.00\% & 0.00\% & 0.00\% & 0.00\% \\
       & 0.37\% & 0.37\% & 0.37\% & 0.00\% & 0.00\% & 0.00\% & 0.00\% & 0.00\% \\
\addlinespace[1ex]
$n=19$ & 28/2/0 & 26/4/0 & 28/2/0 & 28/2/0 & 26/4/0 & 28/2/0 & 29/1/0 & 29/1/0 \\
       & 9.5 & 11.8 & 8.5 & 12.3 & 11.5 & 9.0 & 4.1 & 4.0 \\
       & 2.41\% & 3.97\% & 1.96\% & 2.59\% & 3.74\% & 2.04\% & 1.13\% & 0.75\% \\
       & 0.38\% & 0.26\% & 0.15\% & 0.22\% & 0.14\% & 0.06\% & 0.01\% & 0.01\% \\
\addlinespace[1ex]
$n=24$ & 26/4/0 & 25/5/0 & 26/4/0 & 25/5/0 & 26/4/0 & 26/4/0 & 27/3/0 & 27/3/0 \\
       & 11.9 & 12.3 & 12.3 & 11.5 & 12.7 & 11.7 & 4.9 & 5.3 \\
       & 8.06\% & 8.46\% & 6.49\% & 9.09\% & 8.30\% & 6.34\% & 5.18\% & 5.44\% \\
       & 1.14\% & 0.92\% & 1.23\% & 0.96\% & 0.73\% & 0.62\% & 0.40\% & 0.20\% \\
\midrule
all & 84/6/0 & 81/9/0 & 84/6/0 & 83/7/0 & 82/8/0 & 84/6/0 & 86/4/0 & 86/4/0 \\
    & 7.9 & 8.8 & 7.3 & 8.3 & 8.0 & 6.8 & 3.5 & 3.6 \\
    & 3.49\% & 4.14\% & 2.82\% & 3.89\% & 4.01\% & 2.79\% & 2.11\% & 2.06\% \\
    & 0.63\% & 0.52\% & 0.58\% & 0.39\% & 0.29\% & 0.22\% & 0.14\% & 0.07\% \\
\bottomrule
\end{tabular}
\begin{tablenotes}
\item Number of instances with status Optimal/Feasible/No solution found.
\item Geometric mean of running times (seconds).
\item Average optimality gap, calculated as $100\% \times (\operatorname{UB}^{\operatorname{MIP}}-\operatorname{LB}^{\operatorname{MIP}})/\operatorname{UB}^{\operatorname{MIP}}$.
\item Average gap to optimal, calculated as $100\% \times (\operatorname{UB}^{\operatorname{MIP}}-\operatorname{OPT})/\operatorname{UB}^{\operatorname{MIP}}$.
\end{tablenotes}
\end{threeparttable}
\end{table}

\paragraph{Restricted formulations.}
In three cases, the optimal solution obtained with formulations \textit{Phillips}, \textit{Leung}, and \textit{Leung$^+$} had a cycle time longer than the true optimum due to the modeling restrictions discussed in \Cref{sec:shcsp:pu:remark}.
As a result, the reported gap to the optimum is non-zero, even though the solver status is Optimal and the optimality gap is zero.

\paragraph{Strengthened constraints.}
Interestingly, the advantage of stronger soaking constraints did not persist in the single-thread setting; in fact, in several cases the effect was reversed.
Among the restricted formulations, \textit{Leung} solved three fewer instances to optimality than \textit{Phillips}, which also resulted in an increase in the average running time. 
At the same time, the solutions obtained by \textit{Leung} were often closer to the true optimum, as reflected by the smaller gap to optimality.
For the unrestricted extended formulations, \textit{Imp1} solved one fewer instance to optimality than \textit{Zhou}, and among the base formulations, \textit{Imp2} and \textit{Liu} solved the same number of instances to optimality.
In both comparisons, however, the use of stronger constraints consistently reduced the average gap to the optimum.

\paragraph{Valid inequalities.}
The use of valid inequalities clearly facilitated the solution process.
Specifically, for the \textit{Leung$^+$} formulation, three additional instances were solved to optimality, while \textit{Imp1$^+$} benefited by two instances.
In both cases, the computational time decreased and the gaps to the true optimum were reduced.

\paragraph{Base formulations vs. Extended formulations.}
The base formulations consistently outperformed the corresponding extended versions.
Both for the original and the strengthened soaking constraints, the base formulations (\textit{Liu} and \textit{Imp2}, respectively) solved two more instances to optimality than any of the corresponding extended formulations (\textit{Phillips} and \textit{Zhou}, \textit{Imp1$^+$}, respectively), while requiring significantly less computation time.

\subsubsection{Multiple threads tests}\label{sec:comp:mip8}
\Cref{tab:shcsp:gurobi8} summarizes the results obtained with the Gurobi solver using 8~threads (\texttt{Threads = 8}).

\begin{table}
\centering
\begin{threeparttable}[b]
\caption{Average computational results on randomly generated instances (Gurobi, 8 threads).}
\label{tab:shcsp:gurobi8}
\begin{tabular}{l rrrrrrrr}
\toprule
Instances
 & \multicolumn{1}{c}{Phillips}
 & \multicolumn{1}{c}{Leung}
 & \multicolumn{1}{c}{Leung$^+$}
 & \multicolumn{1}{c}{Zhou}
 & \multicolumn{1}{c}{Imp1}
 & \multicolumn{1}{c}{Imp1$^+$}
 & \multicolumn{1}{c}{Liu}
 & \multicolumn{1}{c}{Imp2}\\
\midrule
$\mu=1.5$ & 30/0/0 & 30/0/0 & 30/0/0 & 30/0/0 & 30/0/0 & 30/0/0 & 30/0/0 & 30/0/0 \\
          & 0.4 & 0.4 & 0.4 & 0.4 & 0.4 & 0.4 & 0.2 & 0.2 \\
          & 0.00\% & 0.00\% & 0.00\% & 0.00\% & 0.00\% & 0.00\% & 0.00\% & 0.00\% \\
          & 0.75\% & 0.75\% & 0.75\% & 0.00\% & 0.00\% & 0.00\% & 0.00\% & 0.00\% \\
\addlinespace[1ex]
$\mu=2.0$ & 30/0/0 & 30/0/0 & 30/0/0 & 30/0/0 & 30/0/0 & 30/0/0 & 30/0/0 & 30/0/0 \\
          & 3.0 & 3.4 & 3.1 & 3.6 & 2.9 & 3.3 & 1.9 & 1.7 \\
          & 0.00\% & 0.00\% & 0.00\% & 0.00\% & 0.00\% & 0.00\% & 0.00\% & 0.00\% \\
          & 0.00\% & 0.00\% & 0.00\% & 0.00\% & 0.00\% & 0.00\% & 0.00\% & 0.00\% \\
\addlinespace[1ex]
$\mu=2.5$ & 27/3/0 & 28/2/0 & 28/2/0 & 27/3/0 & 28/2/0 & 29/1/0 & 29/1/0 & 30/0/0 \\
          & 27.8 & 26.6 & 19.8 & 24.6 & 28.1 & 17.7 & 12.0 & 11.5 \\
          & 3.78\% & 3.16\% & 2.02\% & 3.31\% & 3.34\% & 1.59\% & 0.92\% & 0.00\% \\
          & 0.40\% & 0.27\% & 0.04\% & 0.01\% & 0.04\% & 0.01\% & 0.00\% & 0.00\% \\
\midrule
$n=14$ & 30/0/0 & 30/0/0 & 30/0/0 & 30/0/0 & 30/0/0 & 30/0/0 & 30/0/0 & 30/0/0 \\
       & 1.5 & 1.5 & 1.3 & 1.3 & 1.3 & 1.3 & 0.8 & 0.9 \\
       & 0.00\% & 0.00\% & 0.00\% & 0.00\% & 0.00\% & 0.00\% & 0.00\% & 0.00\% \\
       & 0.37\% & 0.37\% & 0.37\% & 0.00\% & 0.00\% & 0.00\% & 0.00\% & 0.00\% \\
\addlinespace[1ex]
$n=19$ & 30/0/0 & 30/0/0 & 30/0/0 & 30/0/0 & 30/0/0 & 30/0/0 & 30/0/0 & 30/0/0 \\
       & 3.9 & 4.3 & 3.5 & 3.8 & 3.9 & 3.3 & 1.9 & 1.8 \\
       & 0.00\% & 0.00\% & 0.00\% & 0.00\% & 0.00\% & 0.00\% & 0.00\% & 0.00\% \\
       & 0.12\% & 0.12\% & 0.12\% & 0.00\% & 0.00\% & 0.00\% & 0.00\% & 0.00\% \\
\addlinespace[1ex]
$n=24$ & 27/3/0 & 28/2/0 & 28/2/0 & 27/3/0 & 28/2/0 & 29/1/0 & 29/1/0 & 30/0/0 \\
       & 5.1 & 5.4 & 5.3 & 6.4 & 5.8 & 5.9 & 2.6 & 2.7 \\
       & 3.78\% & 3.16\% & 2.02\% & 3.31\% & 3.34\% & 1.59\% & 0.92\% & 0.00\% \\
       & 0.66\% & 0.53\% & 0.30\% & 0.01\% & 0.04\% & 0.01\% & 0.00\% & 0.00\% \\
\midrule
all & 87/3/0 & 88/2/0 & 88/2/0 & 87/3/0 & 88/2/0 & 89/1/0 & 89/1/0 & 90/0/0 \\
    & 3.1 & 3.3 & 2.9 & 3.2 & 3.0 & 3.0 & 1.6 & 1.6 \\
    & 1.26\% & 1.05\% & 0.67\% & 1.10\% & 1.11\% & 0.53\% & 0.31\% & 0.00\% \\
    & 0.38\% & 0.34\% & 0.26\% & 0.00\% & 0.01\% & 0.00\% & 0.00\% & 0.00\% \\
\bottomrule
\end{tabular}
\begin{tablenotes}
\item Number of instances with status Optimal/Feasible/No solution found.
\item Geometric mean of running times (seconds).
\item Average optimality gap, calculated as $100\% \times (\operatorname{UB}^{\operatorname{MIP}}-\operatorname{LB}^{\operatorname{MIP}})/\operatorname{UB}^{\operatorname{MIP}}$.
\item Average gap to optimal, calculated as $100\% \times (\operatorname{UB}^{\operatorname{MIP}}-\operatorname{OPT})/\operatorname{UB}^{\operatorname{MIP}}$.
\end{tablenotes}
\end{threeparttable}
\end{table}

\paragraph{Restricted formulations.}
Similarly to the single-thread experiments, for three cases, the optimal solution obtained with \textit{Phillips}, \textit{Leung}, and \textit{Leung$^+$} had a cycle time longer than the true optimum due to modeling restrictions.

\paragraph{Strengthened constraints.}
The advantage of stronger soaking constraints was visible in the multi-thread setting.
Among the restricted formulations, \textit{Leung} solved one more instance to optimality than \textit{Phillips}, although at the cost of a slightly higher running time.
For the unrestricted extended formulations, \textit{Imp1} solved one more instance to optimality than \textit{Zhou}.
Among the base formulations, \textit{Imp2} solved all instances to optimality.
Note that \textit{Liu} also found the optimal solution for all instances, but failed to prove optimality within the time limit in one case, where the remaining optimality gap was still 27.48\% after ten minutes.

\paragraph{Valid inequalities.}
Although the use of valid inequalities in \textit{Leung$^+$} did not increase the number of instances solved to optimality, it reduced the average running time, the average optimality gap, and the average gap to optimality.
A similar but slightly stronger effect can be observed for \textit{Imp1$^+$}.
In this case, the addition of valid inequalities allowed one more instance to be solved to optimality, while maintaining essentially the same average running time, and reducing both the optimality gaps and the gaps to optimality.

\paragraph{Base formulations vs. extended formulations.}
The base formulations proved to be more effective than the extended versions.
Considering the original soaking constraints, \textit{Liu} outperformed both \textit{Phillips} and \textit{Zhou} in all respects.
Among the formulations using stronger soaking constraints, \textit{Imp2} outperformed \textit{Leung}, 
\textit{Leung$^+$}, \textit{Imp1}, and \textit{Imp1$^+$}.

\subsection{Conclusions}
Overall, the comparison of the MIP~formulations reveals clear differences in formulation strength.
The strengthened formulations generally provide tighter bounds and improved computational performance, while the base formulations consistently outperform the corresponding extended versions.
These results highlight the importance of formulation design, as even structurally similar models can differ substantially in both bound quality and solution efficiency.

\section{Public library for cyclic hoist scheduling}\label{sec:library}
We provide a publicly available library for cyclic hoist scheduling at \url{https://github.com/hmarko89/CyclicHoistScheduling}, with the aim of supporting reproducible research.
The library currently includes:
\begin{enumerate}[(i)]
\item Single-hoist, simple cycle: a natural CP formulation and the MIP formulations of \citet{phillips1976mathematical, liu2002cyclic, zhou2003single, leung2004optimal}, as well as our improved formulations.
\item Single-hoist, multiple cycles: a natural CP formulation and the MIP formulations of \citet{zhou2012mixed, li2014mixed}.
\item Widely used benchmark instances, example instances, and instance generation procedures from the literature.
\end{enumerate}

All MIP formulations are implemented using OR-Tools MathOpt, which allows the instances to be solved with different MIP solvers (e.g., Gurobi, SCIP, HiGHS, GLPK).
Moreover, by exporting the models as LP files, they can be solved with arbitrary MIP solvers.
The CP formulations are implemented using OR-Tools CP-SAT.

Future extensions will include additional MIP and CP formulations, custom solution approaches, and variants for the multi-hoist case, further enhancing the utility of the library for comparative and reproducible research in cyclic hoist scheduling.

\section{Conclusion, future work}
In this paper, we proposed a unified and structured modeling framework for single-hoist cyclic hoist scheduling problems.
The framework establishes a consistent notation and modeling approach that enables a clear and systematic treatment of the problem and its extensions.
Compared to existing formulations in the literature, our approach provides a more transparent and coherent representation of the core constraints, while also supporting a broader range of features, including multifunction tanks, multi-degree cycles, and multi-tank operations.

Building on this framework, we introduced straightforward CP formulations that capture the essential structure of the problem in a natural and reproducible way.
These models can serve as baseline formulations for both simple and extended problem variants.

We revisited MIP formulations from the literature and derived improved variants that integrate their main advantages.
Our computational study demonstrated the effectiveness of the proposed formulations, highlighting the impact of modeling choices on both relaxation strength and solver performance.

To support future research, we provided a publicly available library containing benchmark instances and implementations of several CP and MIP formulations.
This resource, together with the unified modeling framework, is intended to serve as a common reference point for future studies and as a baseline for the development and evaluation of new solution approaches.

Ongoing work focuses on extending the proposed modeling framework to the multi-hoist cyclic scheduling problem, where additional challenges arise from the interaction and coordination of multiple hoists.
Furthermore, we plan to expand the public library with additional formulations, benchmark instances, and solution methods, including both exact and heuristic approaches.

\section*{Acknowledgment}
The author is grateful for the possibility to use HUN-REN Cloud (see \citep{heder2022past}; https://science-cloud.hu/) which helped achieve the results published in this paper.

\section*{Declaration of generative AI and AI-assisted technologies in the manuscript preparation process}
During the preparation of this work, the author used \textit{ChatGPT (OpenAI)} in order to improve language clarity and phrasing, and \textit{Scholar Labs} in order to explore related literature and identify relevant references.
After using these tools, the author reviewed and edited the content as needed and take full responsibility for the content of the published article.

\bibliography{shcsp_references}

\end{document}